\documentclass{article}

\usepackage{amssymb}
\usepackage{hyperref}
\newtheorem{theorem}{Theorem}
\newtheorem{lemma}{Lemma}
\newtheorem{definition}{Definition}

\begin{document}

\title{An $L$-function free proof of Hua's Theorem on sums of five prime squares}\author{CLAUS BAUER\footnote{Claus Bauer, Dolby Laboratories, Beijing, China, cb@dolby.com}}

\date{}

\maketitle

\begin{abstract}
We provide a new proof of Hua's result that every sufficiently large integer $N\equiv 5(mod\,24)$ can be written as the sum of the five prime squares. Hua's original proof relies on the circle method and uses results
from the theory of $L$-functions. Here, we present a proof based on the transference principle first introduced in\cite{green1}. Using a sieve theoretic approach similar to (\cite{shao}), we do not require any results related to the distributions of zeros of $L$- functions. The main technical difficulty of our approach lies in proving the pseudorandomness of the majorant of the characteristic function of the $W$-tricked primes which requires a precise evaluation of the occurring Gaussian sums and Jacobi symbols. \end{abstract}

\section{Introduction}\setcounter{equation}{0}\setcounter{theorem}{0}\setcounter{lemma}{0}
In 1938, Hua (\cite{hua}) showed the following result in additive prime number theory related to the sum of five prime squares:
\begin{theorem}\label{th:hua}
Every sufficiently large integer $N\equiv 5(mod\,24)$ can be written as the sum of five prime squares.
\end{theorem}
Similar to most results in additive prime theory, Hua's proof uses the circle method and relies on the theory of Dirichlet $L$- functions.
Starting with the work of Green (\cite{green1}, \cite{green2}), a different approach to problems in additive prime number theory hat relies on the transference principle in additive combinatorics has been applied. In \cite{li} and \cite{shao1}, it has been used to prove a density version of Vinogradov's three primes theorem. Whereas these proofs rely on results related to the distribution of zeros of $L$-functions, in \cite{shao} a sieve theoretic approach not relying on the theory of $L$-functions is used in combination with the transference principle to derive a new proof of Vinogradov's theorem. The approach in \cite{shao} further differs from the methods in \cite{li},\cite{shao1} by the fact that the transference principle is applied to the set of positive integers $\mathbb{Z}$ instead of being applied to the modular group $\mathbb{Z}/N\mathbb{Z}.$\\
So far, the transference principle has been mostly applied to linear problems in additive prime number theory. In \cite{chow}, a first application to a non-linear problems is described. In this paper, we will prove Theorem \ref{th:hua}. Different from \cite{chow}, we will - as in \cite{shao} - not use the theory of $L$-functions.
A main result and a principal ingredient for the proof of Theorem \ref{th:hua} is the following transference principle in $\mathbb{Z}: $
\begin{theorem}\label{th:1} (Transference principle in $\mathbb{Z}$ with majorant equal to one) Let $0<\delta,\kappa<1,\,1/10<\delta_{1}<1$ be given. Let $N$ be a sufficiently large positive integer. Let $N_{1}=N_{2} =\lfloor N/6\rfloor,$ $N_{3}=N_{4}=\lfloor N/2\rfloor,$ and $N_{5}=N.$ For $i=1,2,3,4,5,$ let $a_{i}:[1,N_{i}]\rightarrow [0,1]$ be arbitrary functions.  Let $\alpha_{i}$ be the average of $a_{i}$ for $i=1,..,5.$ Suppose that the functions $a_{i}$ satisfy the following assumptions:\\
 \\
(1)\,\, (Mean condition) $\alpha_{i}\geq \delta_{1}(i=1,2),\,\alpha_{i}\geq \delta (i=3,4,5),\,\,\frac{1}{2}\left(\min\left(1,\alpha_{3}+\alpha_{4}\right)+\alpha_{4}\right)+\alpha_{5}\geq 1+\delta.$\\
\\
(2)\,\,(Regularity condition for $a_{4}$) \,\,The function $a_{4}$ is $(\delta/50, \kappa)-$ regular.\\
\\
Then
\begin{eqnarray*}
\sum\limits_{n_{1}+n_{2}+n_{3}+n_{4}+n_{5}=N\atop n_{i}\in N_{i}}a_{1}(n_{1})a_{2}(n_{2})a_{3}(n_{3})a_ {4}(n_{4})a_{5}(n_{5})\geq cN^{4},
\end{eqnarray*}
where $c=c(\delta,\kappa)>0$ is a constant depending only on $\delta$ and $\kappa.$
\end{theorem}
The regularity condition (2) in Theorem \ref{th:1} is defined as follows. For $y>2,$ let $P(y)$ be the product of all primes up to $y.$
\begin{definition}\label{de:1}
Let $f:[1,N]\rightarrow \mathbb{R}$ be an arbitrary function. The function is said to be $(\beta,\kappa)$-regular if
\begin{eqnarray*}
\sum\limits_{(u,v)\in M}f(u)f(v)\geq \kappa N^{2},
\end{eqnarray*}
where
\begin{eqnarray*}
M=\left\{(u,v):u\leq \beta N,\,v\geq (1-\beta)N,\,\left(v-u,P(\beta^{-1})\right)=1\right\} .
\end{eqnarray*}
\end{definition}
We note that Theorem \ref{th:1} assumes that the functions $a_{i}$ are bounded by a constant majorant equal to one. For many applications (\cite{green1}, \cite{shao1}, \cite{shao}), this requirement is too strict as one
requires a transference principle valid for more general majorant functions. Due to the assumed discrete majorant property of the functions $a_{i},$ one can derive such transference principles valid for more general majorants that are pseudorandom. Similarly in this paper, for the proof of Theorem \ref{th:hua}, we will require the following transference principle valid for majorants that are pseudorandom:
\begin{theorem}\label{th:1a} (Transference principle in $\mathbb{Z}$)
Let $0<\delta,\kappa<1,\,1/10<\delta_{1}<1$ be given. Let $\eta>0$ be sufficiently small and $N$ be a sufficiently large positive integer. Let $N_{1}=N_{2} =\lfloor N/6\rfloor,$ $N_{3}=N_{4}=\lfloor N/2\rfloor,$ and $N_{5}=N.$ For $i=1,2,3,4,5,$ let $a_{i},v_{i}:[1,N_{i}]\rightarrow \mathbb{R}$ be arbitrary functions.  Let $\alpha_{i}$ be the average of $a_{i}$ for $i=1,..,5.$ Suppose that the functions $a_{i}$ and $v_{i}$ satisfy the following requirements:\\
 \\
(1)\,\,(Majorization condition) $0\leq a_{i}(n)\leq v_{i}(n)\,\,for\,\, all\,\,1\leq n\leq N_{i},\,i=1,2,3,4,5.$\\
\\
(2)\,\, (Mean condition) $\alpha_{i}\geq \delta_{1}(i=1,2),\,\alpha_{i}\geq \delta (i=3,4,5),\,\,\frac{1}{2}\left(\min\left(1,\alpha_{3}+\alpha_{4}\right)+\alpha_{4}\right)+\alpha_{5}\geq 1+\delta.$\\
\\
(3)\,\,(Pseudorandom\,\,condition) The functions $v_{i},\,i=1,..,5$ \,\,are\,\,pseudorandom.\\
\\
(4)\,\,(Discrete\,\,majorant\,\,property)\,\,The functions $a_{i},\,i=1,..,5$ satisfy the discrete majorant property for some $4<q<5.$ \\
\\
(5)\,\,(Regularity condition for $a_{4}$) \,\,The function $a_{4}$ is $(\delta/50, \kappa)-$ regular.\\
\\
Then
\begin{eqnarray*}
\sum\limits_{n_{1}+n_{2}+n_{3}+n_{4}+n_{5}=N\atop n_{i}\in N_{i}}a_{1}(n_{1})a_{2}(n_{2})a_{3}(n_{3})a_ {4}(n_{4})a_{5}(n_{5})\geq cN^{4},
\end{eqnarray*}
where $c=c(\delta,\kappa)>0$ is a constant depending only on $\delta$ and $\kappa.$
\end{theorem}
To define the pseudorandom condition (3) and the discrete majorant property (4) in Theorem \ref{th:1a}, we introduce some more terminology:
For a (compactly supported) function $f:\mathbb{Z}\rightarrow \mathbb{R},$ its Fourier transform is defined as
\begin{eqnarray*}
\hat{f}(\theta)=\sum\limits_{n\in \mathbb{Z}}f(n)e(n\theta),
\end{eqnarray*}
where $e(n\theta)=exp(2\pi i n\theta).$ The $L^{q}$ norm of the Fourier transform is defined by
\begin{eqnarray*}
||\hat{f}||_{q} =\left(\int\limits_{0}^{1}|\hat{f}(\theta)|^{q}d\theta\right)^{1/q}.
\end{eqnarray*}
\begin{definition}
The function $f$ is said to be $\eta$-pseudoranndom  if $|\hat{f}(r/N)-\delta_{r,0}N|\ll \eta N$ for each $r\in \mathbb{Z}/N\mathbb{\mathbb{Z}},$ where $\delta_{r,0}$ is the Kronecker delta.
\end{definition}
\begin{definition}
The function $f$ is said to satisfy the discrete majorant property if $||\hat{f}||_{q}\ll_{q}N^{1-1/q},$ where the implied constant depends on $q$ only.
\end{definition}
Our paper is structured as follows. We first prove Theorem \ref{th:1} in section \ref{se:th1}. Subsequently, we derive Theorem \ref{th:1a} from Theorem \ref{th:1} in section \ref{se:th1a}. In section \ref{se:def}, we will define the {\it $W$-tricked primes squares} and their majorant function. In sections \ref{se:mean} - \ref{se:regularity}, we will show that the characteristic function of the {\it $W$-tricked primes squares} and its majorant satisfy the assumptions of Theorem \ref{th:1a}. This will allow us applying Theorem \ref{th:1a} to prove Theorem \ref{th:hua}.
\section{Proof of Theorem \ref{th:1}}\label{se:th1}
\setcounter{equation}{0}\setcounter{theorem}{0}\setcounter{lemma}{0}
We will make use of the following Theorem (\cite[Th. 1.4]{shao}):
\begin{theorem}\label{th:2} (Transference principle in $\mathbb{Z}$ with majorant equal to one) Let $0<\delta,\kappa<1$ be given. Then for sufficiently small $\eta >0$ and sufficiently large positive $N,$ the following statement holds: Let $N_{1}=N_{2}=\lfloor N/2\rfloor ,$ and $N_{3}=N.$ For $i=1,2,3$ let $a_{i}:[1,N_{i}]\rightarrow [0,1]$ be arbitrary functions.
 Let $\alpha_{i}$ be the average of $a_{i}$ for $i=1,2,3.$ Suppose that the functions $a_{i}$ satisfy the following assumptions:\\
 \\
(1)\,\, (Mean condition) $\alpha_{i}\geq \delta (i=1,2,3),\,\,\frac{1}{2}\left(\min\left(1,\alpha_{1}+\alpha_{2}\right)+\alpha_{2}\right)+\alpha_{3}\geq 1+\delta.$\\
\\
(2)\,\,(Regularity condition for $a_{3}$) \,\,The function $a_{1}$ is $(\delta/50, \kappa)-$ regular.\\
\\
Then
\begin{eqnarray*}
\sum\limits_{n_{1}+n_{2}+n_{3}=N\atop n_{i}\in N_{i}}a_{1}(n_{1})a_{2}(n_{2})a_{3}(n_{3})\geq cN^{2},
\end{eqnarray*}
where $c=c(\delta,\kappa)>0$ is a constant depending only on $\delta$ and $\kappa.$\end{theorem}
We now proceed to the proof of Theorem \ref{th:1}. For a given positive integer $m\leq N$ we define $m_{i},\,i=1,2,3,4,5,$ as a function of $m$ in the same way we have defined $N_{i}$ as a function of $N$ in Theorem \ref{th:1}.
Further, we for each $i=1,2,3,4,5,$ we define a function $a_{i,m}(n)$ as follows: $a_{i,m}(n)=a_{i}(n)$ for $n\leq m_{i}$ and $a_{i,m}(n)=0$ for $m_{i}<n\leq N_{i}.$ Thus $a_{i,m}(n)\leq a_{i}(n),$ $\forall \,n\in [1,N_{i}].$\\
We note that Theorem \ref{th:2} is valid for  all $N>N_{c},$ where $N_{c}$ is a large positive integer. If we choose $N$ sufficiently large, we can assume that $N/100>N_{c},$ which we will do in the following. Applying
Theorem \ref{th:2} to the functions $a_{3,m}, a_{4,m}$ and $a_{5,m}$ and using the mean condition $\alpha_{i}\geq \frac{1}{10},\,i=1,2,$ we see
\begin{eqnarray*}
&&\sum\limits_{n_{1}+n_{2}+n_{3}+n_{4}+n_{5}=N\atop n_{i}\leq N_{i}}a_{1}(n_{1})a_{2}(n_{2})a_{3}(n_{3})a_{4}(n_{4})a_{5}(n_{5})\nonumber\\
&=&\sum\limits_{m\leq N}\sum\limits_{n_{1}+n_{2}=N-m\atop n_{i}\leq N_{i}}a_{1}(n_{1})a_{2}(n_{2})
\sum\limits_{n_{3}+n_{4}+n_{5}=m\atop n_{i}\leq N_{i}}a_{3}(n_{3})a_{4}(n_{4})a_{5}(n_{5})\nonumber\\
&\geq &\sum\limits_{N/100<m\leq N}\sum\limits_{n_{1}+n_{2}=N-m\atop n_{i}\leq N_{i}}a_{1}(n_{1})a_{2}(n_{2})
\sum\limits_{n_{3}+n_{4}+n_{5}=m\atop n_{i}\leq m_{i}}a_{3,m}(n_{3})a_{4,m}(n_{4})a_{5,m}(n_{5})\nonumber\\
&\geq & \frac{cN^{2}}{10000}\sum\limits_{N/100<m\leq N}\sum\limits_{n_{1}+n_{2}=N-m\atop n_{i}\leq N_{i}}a_{1}(n_{1})a_{2}(n_{2})\\
&\geq & \frac{cN^{2}}{10000}\sum\limits_{2N/3<m\leq 99N/100}\sum\limits_{n_{1}+n_{2}=N-m\atop n_{i}\in [N/200,N_{i}]}a_{1}(n_{1})a_{2}(n_{2})\\
&=& \frac{cN^{2}}{10000}\sum\limits_{ n_{1},n_{2}\in [N/200,N_{i}]}a_{1}(n_{1})a_{2}(n_{2})\\
&\geq & \frac{cN^{4}}{10000}\left(\delta_{1}-\frac{1}{200}\right)^{2},
\end{eqnarray*}
which proves Theorem \ref{th:1} with $c:=10^{-8}c.$
\section{Proof of Theorem \ref{th:1a}}\label{se:th1a}
\setcounter{equation}{0}\setcounter{theorem}{0}\setcounter{lemma}{0}
We introduce some additional terminology. For notational convenience, we will fix some $i=1,2,3,4,5,$ and set $a(n):=a_{i}(n),\,\alpha_{i}:=\alpha,\,N:=N_{i}.$  Let $0<\epsilon<1$ be a small parameter chosen later which depends only on $\delta$ and $\kappa.$ We set $\mathbb{T}:=[0,1].$ Let
\begin{eqnarray*}
T_{\epsilon}=\{\theta\in \mathbb{T}:|\hat{a}(n)|\geq \epsilon N\}.
\end{eqnarray*} Define
\begin{eqnarray*}
B=B_{\epsilon}=\{1\leq b\leq \epsilon N:||b\theta||<\epsilon\,\,for\,\,all\,\,\theta\in T_{\epsilon}\},
\end{eqnarray*} where $||x||$ denotes the distance from $x$ to its closest integer.
Using these definitions, we define the following functions:
\begin{eqnarray*}
a'(n)=\mathbb{E}_{b_{1},b_{2}\in B}a(n+b_{1}-b_{2})=\frac{1}{B^{2}}\sum\limits_{b_{1},b_{2}\in B}a(n+b_{1}-b_{2}),\,\,a''(n)=a(n)-a'(n).
\end{eqnarray*} We will need the following Lemma:
\begin{lemma}\label{le:aa}Suppose that $\eta$ is sufficiently small depending on $\epsilon.$ The functions $a'$ and $a''$ defined above have the following properties:\\
\\
(1)\, $(a'\,\, is\,\, set-like)\,\, 0\leq a'(n)\leq 1+O_{\epsilon}(\eta)\,\,for\,\,any\,\,n.$ Moreover, $\mathbb{E}_{1\leq n\leq N}a'(n)=\alpha+O(\epsilon).$\\
\\
(2)\, $(a''\,\,is\,\,uniform)\,\,a''(\theta)=O(\epsilon N)\,\,for\,\,all\,\,\theta.$\\
\\
(3)$\,(a_{4}'\,is\,regular)\,a'_{4}\,is\,\,(\delta/50,\kappa-O(\epsilon))-regular.$\\
\\
(4)$\,\,||\hat{a}'||_{q}\leq ||\hat{a}||_{q},\,\,||\hat{a}''||_{q}\leq ||\hat{a}||_{q}.$
\end{lemma}
{\it Proof:} The proof is word-by-word identical with the proof of \cite[Lemma 3.3]{shao}.\\
\\ The following Lemma is the main ingredient for the proof of Theorem \ref{th:1a}:
\begin{lemma}\label{le:aaa}
With the functions $a_{i}$ and $a'_{i}$ as defined above, we have
\begin{eqnarray*}
&&\Big| \sum\limits_{n_{i}\in N_{i},\,i=1,,,,5\atop n_{1}+n_{2}+n_{3}+n_{4}+n_{5}=N}a_{1}(n_{1})a_{2}(n_{2})a_{3}(n_{3})a_{4}(n_{4})a_{5}(n_{5})\nonumber\\
 &&- \sum\limits_{n_{i}\in N_{i},\,i=1,,,,5\atop n_{1}+n_{2}+n_{3}+n_{4}+n_{5}=N}a'_{1}(n_{1})a'_{2}(n_{2})a'_{3}(n_{3})a'_{4}(n_{4})a'_{5}(n_{5})\Big|\ll \epsilon^{5-q}N^{4}.
\end{eqnarray*}
\end{lemma}
Proof: The difference on the left side can be expressed as the sum of several terms of the form
\begin{eqnarray}\label{eq:sec71mm}
\sum\limits_{n_{i}\in N_{i},\,i=1,,,,5\atop n_{1}+n_{2}+n_{3}+n_{4}+n_{5}=N}f_{1}(n_{1})f_{2}(n_{2})f_{3}(n_{3})f_{4}(n_{4})f_{5}(n_{5})=\int\limits_{0}^{1}\hat{f_{1}}(\theta)
\hat{f_{2}}(\theta)\hat{f_{3}}(\theta)\hat{f_{4}}(\theta)\hat{f_{5}}(\theta)e(-N\theta)\,d\theta,
\end{eqnarray}
where
$f_{i}\in\{a_{i},a'_{i},a''_{i}\},$ and $f_{i}=a''_{i}$ for at least one $i.$ Without loss of generality, we assume that $f_{5}=a''_{5}.$ By H\"{o}lder's inequality, for some $4<q<5,$ the right-hand side of (\ref{eq:sec71mm}) is bounded by
\begin{eqnarray*}
||\hat{f_{5}}||_{\infty}^{5-q}||\hat{f_{5}}||_{q}^{q-4}||\hat{f_{1}}||_{q}||\hat{f_{2}}||_{q}||\hat{f_{3}}||_{q}||\hat{f_{4}}||_{q}.
\end{eqnarray*} By Lemma \ref{le:aa}, $||\hat{f}_{5}||_{\infty}\ll \epsilon N.$ Further, by the discrete majorant property and Lemma \ref{le:aa}, $||\hat{f}_{i}||_{q}\ll N^{1-\frac{1}{q}},\,i=1,2,3,4,5.$ Combining these estimates, we obtain the desired estimate.\\
\\
\\
{\it Proof of Theorem \ref{th:1a}:}
By Lemma \ref{le:aa}, the functions $a'_{i}$ are all bounded above uniformly by $1+O_{\epsilon}(\eta)$ with averages $\alpha_{i}+O(\epsilon)$ and $a'_{4}$ is $(\delta/50,\kappa/2)$ - regular.
If $\epsilon$ and $\eta$ are chosen small enough, then Theorem \ref{th:1} implies
\begin{eqnarray}\label{eq:b2mm}
\sum\limits_{n_{1}+n_{2}+n_{3}+n_{4}+n_{5}=N\atop n_{i}\in N_{i}}
a'_{1}(n_{1})a'_{2}(n_{2})a'_{3}(n_{3})a'_{4}(n_{4})a'_{5}(n_{5}) \geq c'N^{4}.
\end{eqnarray}
Combining this with Lemma \ref{le:aaa}, we deduce by choosing $\epsilon$ small enough that
\begin{eqnarray*}
\sum\limits_{n_{1}+n_{2}+n_{3}+n_{4}+n_{5}=N\atop n_{i}\in N_{i}}
a_{1}(n_{1})a_{2}(n_{2})a_{3}(n_{3})a_{4}(n_{4})a_{5}(n_{5}) \geq \frac{c'}{2}N^{4}.
\end{eqnarray*}
This completes the proof of Theorem \ref{th:1a}.
\section{Proof of Theorem \ref{th:hua}}\label{se:def}
We first introduce the concept of $W$- tricked prime squares and define their majorant function. Subsequently, we prove Theorem \ref{th:hua}.
\subsection{$W$-tricked prime squares and Selberg's majorant}\label{se:char}
\setcounter{equation}{0}\setcounter{theorem}{0}\setcounter{lemma}{0}
Let $W=8\prod\limits_{2<p\leq w}p$ be the product of the number eight and all odd primes not larger than a constant $w$ and let $b$ be a reduced residue class modulo $W.$ We set $[W]:=\{1,..,W\}.$ We say that an integer $q$ is $w$-smooth if none of its prime divisors is larger than $w.$ Otherwise, we call $q$ $w$-rough. We denote the set of primes by $\mathbb{P}$ and the set of prime squares by $\mathbb{P}^{2}.$ We fix a small constant $\delta$ and set $z_{i}=N_{i}^{1/4-2\delta},$ $i=1,...,5,$ where $N_{i}$ is defined as in Theorem \ref{th:1}. We assume that $N$ is a sufficiently large integer depending on $w$ and $\delta.$
Further, we define $P:=P_{i}$ as the product of all primes $p<z_{i}$ and $(p,W)=1.$ We set
\begin{eqnarray*}
\sigma(b)&:=\sharp&\left\{z\in [W]:z^{2}-b\equiv 0(mod\,W)\right\}.
\end{eqnarray*} We define the $\sigma(b)$ reduced residue classes $h_{j},1\leq j\leq \sigma(b),$  modulo $W$ via the relation $h_{j}^{2}\equiv b(mod\,W).$
Using ideas from \cite{browning} and \cite{shao}, we define the $W$-shifted prime squares:
\begin{eqnarray*}
a_{i}(n)&=&\left\{\begin{array}{ll}
\frac{2\phi(W)\sqrt{Wn+b_{i}}\log\,z_{i}}{W\sigma(b)},& if\,\,Wn+b_{i}=x^{2}\,for\,some\,prime\,number\\ &x \,with\,z_{i}\leq x\leq \sqrt{WN_{i}+b_{i}},\\0,& else.
\end{array}\right\}\end{eqnarray*}
Following the approach in  \cite{shao}, we use Selberg's upper bound sieve to define the majorant of $a_{i}(n)$ as follows:
\begin{eqnarray}\label{eq:v}
v_{i}(n)&=&\left\{\begin{array}{ll}
\frac{2\phi(W)\sqrt{Wn+b_{i}}\log\,z_{i}}{W\sigma(b)}\left(\sum\limits_{d|(\sqrt{Wn+b_{i}},P)}p_{d}\right)^{2},& if\,\,Wn+b_{i}=x^{2}\,for\\&\,some\,natural\,
integer\,x\leq \sqrt{WN_{i}+b_{i}},\\0,& else.
\end{array}\right\}\nonumber\\
\end{eqnarray}
The real weights $\rho_{d}$ are supported on $d<z_{i},\,\mu(d)\ne 0,$ and satisfy $|\rho_{d}|\leq 1,$ and $\rho_{1}=1.$
If $a_{i}(n)\ne 0,$ then $\sqrt{Wn+b_{i}}\geq z_{i}$ and therefore $(Wn+b_{i},P)=1,$ such that $a_{i}(n)=v_{i}(n).$ This shows that $v_{i}(n)$ is a majorant of $a_{i}(n).$\\
\\
For later usage, we introduce the following notation:
\begin{eqnarray}\label{eq:defj}
J:=J_{i}=\sum\limits_{d|P\atop d<z_{i}}\frac{1}{\phi(d)} =\sum\limits_{d<z_{i},\,d\,squarefree\atop (d,W)=1}\frac{1}{\phi(d)}.
\end{eqnarray}
We know from \cite[Appendix A, (A.1)]{shao} that
\begin{eqnarray}\label{eq:kjjj}
\label{eq:J} J_{i}=\frac{\phi(W)}{W}(\log\,z_{i}+O_{W}(1)).
\end{eqnarray}
\subsection{Proof of Theorem \ref{th:hua}}
We will need the following auxiliary Lemmas:
\begin{lemma}\label{le:cdc} Let $p$ be a prime number.
Let $A,\,B,\,C$ be non-empty subsets of $\mathbb{Z}/p\mathbb{Z}$ such that $|A|+|B|+|C|\geq p+2.$ Then,
\begin{eqnarray*}
A+B+C=\mathbb{Z}/p\mathbb{Z}.
\end{eqnarray*}
\end{lemma}
{\it Proof:} See \cite[Proposition 2.6]{tao}.
\begin{lemma}\label{le:pp}
Let $p$ be a prime number $\geq 5$ and define $A_{p}$ as the set of all quadratic residues modulo $p.$
Then ,
\begin{eqnarray*}
A_{p}+A_{p}+A_{p}+A_{p}+A_{p}=\mathbb{Z}/p\mathbb{Z}.
\end{eqnarray*}
\end{lemma}
{\it Proof of Lemma \ref{le:pp}:} As $|A_{p}|=(p-1)/2,$ we see that for $p\geq 7,$
\begin{eqnarray*}
|A_{p}+A_{p}+A_{p}|+|A_{p}|+|A_{p}|
\geq |A_{p}|+|A_{p}|+|A_{p}|\geq p+2.
\end{eqnarray*}
For $p\geq 7,$ the Lemma now follows from Lemma \ref{le:cdc}. For $p=5,$ the Lemma follows by a case by case analysis.
\begin{lemma}\label{le:MM}
For any integer $N\equiv 5(mod\,24),$ there are integers $b_{i},\,1\leq b_{i}\leq W,\,(W,b_{i})=1,\,i=1,...,5$ where each $b_{i}$ is a quadratic residue modulo $W$ such that
$N\equiv b_{1}+b_{2}+b_{3}+b_{4}+b_{5}(mod\,W).$
\end{lemma}
{\it Proof:} By case-by-case inspection, we see that $b=1$ is the only quadratic residue modulo $24.$ Therefore, $N\equiv b_{1}+b_{2}+b_{3}+b_{4}+b_{5}(mod\,24)$ for any quadratic residues $b_{i}$ modulo $W.$ The Lemma now follows from Lemma \ref{le:pp} and the Chinese remainder theorem.\\
\\
{\it Proof of Theorem \ref{th:hua}}:
Let $M$ be a sufficiently large integer with $M\equiv 5(mod\,24).$ We will show that $M$ can be written as the sum of five prime squares. By Lemma \ref{le:MM}, we can find integers $b_{1},...,b_{5}$ which are quadratic residues modulo $W$ such that $M\equiv b_{1}+b_{2}+b_{3}+b_{4}+b_{5}(mod\,W).$
 We set $N=(M-b_{1}-b_{2}-b_{3}-b_{4}-b_{5})/W.$ Using the integers $b_{i},$
we define the functions $a_{i}$ and $v_{i}$ as in section \ref{se:char}.
 In section \ref{se:char}, we have shown that the functions $v_{i}$ are majorants for the functions $a_{i}.$
In sections \ref{se:mean} - \ref{se:regularity},
we will show that the functions $a_{i}$ and $v_{i}$ also satisfy the conditions (2) - (5) of Theorem \ref{th:1a}. Therefore, by Theorem \ref{th:1a} there exist $n_{i}\in N_{i},\,i=1,...,5,$ such that $a_{i}(n_{i})>0,$ i.e., $Wn_{i}+b_{i}\in \mathbb{P}^{2},$ and
$N=n_{1}+n_{2}+n_{3}+n_{4}+n_{5}.$ Thus,
\begin{eqnarray*}
M=WN+\sum\limits_{i=1}^{5}b_{i}= \sum\limits_{i=1}^{5}Wn_{i}+b_{i},
\end{eqnarray*} which proves that $M$ is the sum of five prime squares.
\section{Mean condition}\label{se:mean}
\begin{lemma}\label{le:meancond} For $i=1,...,5,$ we define by $\alpha_{i}$ the mean of $a_{i}(n).$\\
i)\,\,For $i=1,...,5,$
\begin{eqnarray*}\label{eq:lk}\alpha_{i}\geq 1/2-5\delta.\end{eqnarray*}
ii) For $\delta=0.001,$ the functions $a_{i}(n),\,i=1,...5,$ satisfy the mean condition (2) of Theorem \ref{th:1a}.
\end{lemma}
{\it Proof of Lemma \ref{le:meancond}:} To simplify notation, we write $a(n)=a_{i}(n),$ $b(n)=b_{i}(n),$ $z=z_{i},$ and $N=N_{i}.$
We first prove part i).
\begin{eqnarray*}
&&\sum\limits_{n\leq N}a(n)=\frac{2\phi(W)}{W\sigma(b)}\log\,z\sum\limits_{n\leq N\atop {Wn+b=x^{2} \,for\, some\, prime\, x \,with \atop z\leq x\leq\sqrt{WN+b} }}\sqrt{Wn+b}\nonumber\\
&=& \frac{2\phi(W)}{W\sigma(b)}\log\,z\sum\limits_{z\leq x\leq \sqrt{WN+b},\,x\,prime \atop x^{2}-b\equiv 0(mod\,W)}x=\frac{2\phi(W)}{W\sigma(b)}\log\,z\sum\limits_{j=1}^{\sigma(b)}\sum\limits_{z\leq x\leq \sqrt{WN+b},\,x\,prime\atop x\equiv h_{j}(mod\,W)}x\nonumber\\
&\geq &\frac{2\phi(W)}{W\sigma(b)}\log\,z\sum\limits_{j=1}^{\sigma(b)}\left(\frac{NW}{2\phi(W)\log (\sqrt{WN+b})}1-\epsilon\right)\nonumber\\
&= & N\log\,z/\log (\sqrt{WN+b})-2\epsilon\log\,z\geq N(1/4-2\delta)2(1-\epsilon)\nonumber\\
&\geq &N(1/2-5\delta),
\end{eqnarray*}
for $\epsilon\leq \delta/100$
which implies part i) of the Lemma. Part ii) is a direct consequence of i).
\section{Pseudorandom condition}\setcounter{equation}{0}\setcounter{theorem}{0}\setcounter{lemma}{0}
In this section, we set $v(n)=v_{i}(n),$ $N=N_{i},$ $z=z_{i},$ and $J=J_{i},$ for fixed $i\in\{1,2,3,4,5\}.$ The main purpose of this section is to prove the following Lemma:
\begin{lemma}\label{le:pseudo} For any $r\in \mathbb{Z}/N\mathbb{Z},$
\begin{eqnarray*}
\hat{v}\left(\frac{r}{N}\right)=\left(\delta_{r,0}+O_{w}\left(w^{-1/2+\epsilon}\right)\right)N,
\end{eqnarray*}
where $\delta_{r,0}$ is the Kronecker delta.
\end{lemma}
To prepare for the proof of Lemma \ref{le:pseudo}, we introduce some further terminology. For integers $q,d_{1},d_{2},W,$ we define $q_{d_{1},d_{2}}:=q/(q,[d_{1},d_{2}]^{2})$ and
$a_{d_{1},d_{2}}=a[d_{1},d_{2}]^{2}/(q,[d_{1},d_{2}]^{2}).$ We note that $(W,q)=(W,q_{d_{1},q_{2}})$ is independent of $[d_{1},d_{2}]$ as $(W,d_{1}d_{2})=1$ which follows from the summation condition
$d_{i}|\sqrt{Wn+b}$ in (\ref{eq:v}). We define,
\begin{eqnarray}\label{eq:w3w3w3}
h:=(W,q)=(W,q_{d_{1},q_{2}}),\qquad W_{1}:=W/h,,\qquad q_{W,d_{1},d_{2}}:={q_{d_{1},d_{2}}}/h.
\end{eqnarray}
For the proof of Lemma \ref{le:pseudo}, we will analyse the following term defined for $1\leq Y\leq N:$
\begin{eqnarray}\label{eq:f}
f_{d_{1},d_{2}}(Y,\alpha)&:=&\sum\limits_{n\leq Y\atop{ [d_{1},d_{2}]|\sqrt{Wn+b}\atop Wn+b=x^{2}\,for\,some\,x\leq \sqrt{WY+b}}}
\frac{2\sqrt{Wn+b}}{\sigma(b)}e\left(\alpha n \right).
\end{eqnarray}
For the analysis of $f_{d_{1},d_{2}}(Y,\alpha)$, we divide the integral $[0,W\lfloor z^{2}\rfloor]$ into major arcs $\mathfrak{M}$ and minor arcs $\mathfrak{m}$ as follows: We set $Q=\lfloor{N^{2\delta/5}}\rfloor,$ $R=\lfloor N^{1-\delta/2}\rfloor, $ and
\begin{eqnarray}\label{eq:arcs}
\mathfrak{M}&=&\bigcup_{q_{1}\leq Q}\bigcup_{a_{1}=1\atop (a_{1},q_{1})=1}^{W\lfloor z^{2}\rfloor}\mathfrak{M}(a_{1},q_{1}),\nonumber\\ \mathfrak{M}(a_{1},q_{1})&=&\left\{\alpha\in [0,W\lfloor z^{2}\rfloor]:\left|\alpha-\frac{a_{1}}{q_{1}}\right|\leq \frac{1}{q_{1}R}\right\},\nonumber\\
\mathfrak{m}&=&[0,W\lfloor z^{2}\rfloor]\setminus \mathfrak{M}.
\end{eqnarray} Using the definitions (\ref{eq:f}) and (\ref{eq:arcs}), we can express $\hat{v}\left(\frac{r}{N}\right)$ as follows:
\begin{eqnarray}\label{eq:defv1v2}
\hat{v}\left(\frac{r}{N}\right)=\hat{v}_{1}\left(\frac{r}{N}\right)+\hat{v}_{2}\left(\frac{r}{N}\right)+\hat{v}_{3}\left(\frac{r}{N}\right),
\end{eqnarray}
where
\begin{eqnarray}\label{eq:defv1}
\hat{v}_{1}\left(\frac{r}{N}\right)&=&\frac{\phi(W)\log\,z}{W}\sum\limits_{d_{1},d_{2}|P\atop q| [d_{1},d_{2}]^{2}}p_{d_{1}}p_{d_{2}}f_{d_{1},d_{2}}(N,r/N),\nonumber\\
\hat{v}_{2}\left(\frac{r}{N}\right)&=&\frac{\phi(W)\log\,z}{W}\sum\limits_{d_{1},d_{2}|P\atop{q\nmid [d_{1},d_{2}]^{2}\atop rW[d_{1},d_{2}]^{2}/N\in \mathfrak{M}}}p_{d_{1}}p_{d_{2}}f_{d_{1},d_{2}}(N,r/N),\nonumber\\
\hat{v}_{3}\left(\frac{r}{N}\right)&=&\frac{\phi(W)\log\,z}{W}\sum\limits_{d_{1},d_{2}|P\atop{q\nmid [d_{1},d_{2}]^{2}\atop rW[d_{1},d_{2}]^{2}/N\in \mathfrak{m}}}p_{d_{1}}p_{d_{2}}f_{d_{1},d_{2}}(N,r/N).
\end{eqnarray}
Lemma \ref{le:pseudo} follows from the following three Lemmas:
\begin{lemma}\label{le:pseudo1} For any $r\in \mathbb{Z}/N\mathbb{Z},$
\begin{eqnarray*}
\hat{v}_{1}\left(\frac{r}{N}\right)=\left(\delta_{r,0}+O_{w}\left(w^{-1/2+\epsilon}\right)\right)N.
\end{eqnarray*}
\end{lemma}
\begin{lemma}\label{le:pseudo1a} For any $r\in \mathbb{Z}/N\mathbb{Z},$
\begin{eqnarray*}
\hat{v}_{2}\left(\frac{r}{N}\right)\ll_{w}Nw^{-1/2+\epsilon}.
\end{eqnarray*}
\end{lemma}
\begin{lemma}\label{le:pseudo2} For any $r\in \mathbb{Z}/N\mathbb{Z},$
\begin{eqnarray*}
\hat{v}_{3}\left(\frac{r}{N}\right)\ll N^{1-\delta/400}.
\end{eqnarray*}
\end{lemma}
We prove the Lemmas \ref{le:pseudo1} - \ref{le:pseudo2} in the next three sections.  We will use ideas from \cite[Section 5]{browning} and \cite[Appendix A]{shao}.
\subsection{Proof of Lemma \ref{le:pseudo1}}
We initially assume $\frac{r}{N}=\frac{a}{q}$
and analyse the term of $f_{d_{1},d_{2}}(Y,a/q)$ in section \ref{se:ll}
Then, we prove Lemma \ref{le:pseudo1} in section \ref{se:ko}.
\subsubsection{The case $\frac{r}{N}=\frac{a}{q}$}\label{se:ll}
We will make use of the following lemma:
\begin{lemma}\label{le:a3} For any positive integer $q$ dividing $P,$ the sum
\begin{eqnarray}\label{eq:tq}
T(q):=\sum\limits_{d_{1},d_{2}|P\atop q|[d_{1},d_{2}]}
\frac{p_{d_{1}}p_{d_{2}}}{[d_{1},d_{2}]}
\end{eqnarray}
satisfies
\begin{eqnarray*}
|T(q)|\ll J^{-1}q^{-1+\epsilon},
\end{eqnarray*}
where $J$ is defined in (\ref{eq:defj}).
Moreover,
$T(1)=J^{-1}.
$
\end{lemma} {\it Proof of Lemma:} See \cite[Lemma A.3]{shao}.
\begin{lemma}\label{le:1}
For integers $q,d_{1},d_{2},W,$ with $q|[d_{1},d_{2}]^{2},$ $([d_{1},d_{2}],W)=1,$ and for two co-prime integers $a$ and $q$ there is:
\begin{eqnarray*}
f_{d_{1},d_{2}}(Y,a/q)&=&\frac{\epsilon_{q}Y}{[d_{1},d_{2}]}+O\left(Y^{1/2}W^{1/2}\right),
\end{eqnarray*}
where $\epsilon_{q}=\epsilon_{q}(a/q,W,b)$ does not depend on $Y.$ Moreover $\epsilon_{q}=1$ if $q=1,$ $|\epsilon_{q}|\leq 1$ if $q>1,\,(q,W)=1,$ and $\epsilon_{q}=0$ if $(q,W)>1.$
\end{lemma}
{\it Proof of Lemma \ref{le:1}:}
We write $f_{d_{1},d_{2}}\left(Y,a/q\right)$ as follows:
\begin{eqnarray}\label{eq:pseudo11111}
f_{d_{1},d_{2}}\left(Y,a/q\right)
&=&\frac{2}{\sigma(b)}\sum\limits_{x\leq \sqrt{WY+b}\atop {[d_{1},d_{2}]|x\atop x^{2}-b\equiv 0(mod\,W)}}x
 e\left(\frac{a\left( x^{2}-b\right)}{qW}\right)\\
 \label{eq:pseudo1}  &=&\frac{2e(\frac{-a b}{qW})[d_{1},d_{2}]}{\sigma(b)}\sum\limits_{u\leq \sqrt{WY+b}/[d_{1},d_{2}]\atop u^{2}[d_{1},d_{2}]^{2}- b\equiv 0(mod\,W)}u
 e\left(\frac{a[d_{1},d_{2}]^{2}u^{2}}{qW}\right).\nonumber\\
\end{eqnarray}
Breaking the sum over $u$ into congruence classes modulo $W,$ we see from (\ref{eq:pseudo1}) and $q|[d_{1},d_{2}]^{2}:$
\begin{eqnarray}\label{eq:pseudo1a}
&&f_{d_{1},d_{2}}\left(Y,a/q\right)\nonumber\\
  &=&\frac{2e\left(-\frac{a b}{qW}\right)[d_{1},d_{2}]}{\sigma(b)}\sum\limits_{z\in [W]\atop z^{2}[d_{1},d_{2}]^{2}-b\equiv 0(mod\,W)}\sum\limits_{z+yW\leq \sqrt{WY+b}/[d_{1},d_{2}]}(z+yW)
 e\left(\frac{a[d_{1},d_{2}]^{2} (z+yW)^{2}}{qW}\right)\nonumber\\
   &=&\frac{2e\left(\frac{-a b}{qW}\right)[d_{1},d_{2}]}{\sigma(b)}\sum\limits_{z\in [W]\atop z^{2}[d_{1},d_{2}]^{2}-b\equiv 0(mod\,W)}  e\left(\frac{a[d_{1},d_{2}]^{2}z^{2}}{qW}\right)\sum\limits_{z+yW\leq \sqrt{WY+b}/[d_{1},d_{2}]}(z+yW)\nonumber\\
     &=&\left(\frac{Y}{[d_{1},d_{2}]}+O\left(Y^{1/2}W^{1/2}\right)\right)\frac{1}{\sigma(b)}\sum\limits_{z\in [W]\atop z^{2}[d_{1},d_{2}]^{2}-b\equiv 0(mod\,W)}  e\left(\frac{a([d_{1},d_{2}]^{2}z^{2}-b)}{qW}\right).
\end{eqnarray}
We note that the inner sum over $z$ is equal to $\sigma(b\overline{[d_{1},d_{2}]^{2}})=\sigma(b)$ if $q=1.$
If $q>1,$ and $(q,W)=1$ the absolute value of the sum over $z$ is $\leq \sigma(b\overline{[d_{1},d_{2}]^{2}})=\sigma(b).$
Further, if $(q,W)>1,$ the assumptions $(W,b)=1$ and $q|[d_{1},d_{2}]^{2},$ contradict with the summation condition $[d_{1},d_{2}]|\sqrt{Wn+b}$ in (\ref{eq:f}), i.e.,
$(q,W)>1\Rightarrow  f_{d_{1},d_{2}}\left(Y,a/q\right)=0.$ The Lemma now follows from (\ref{eq:pseudo1a}).\\
\\
By the definition of $v(n)$ in (\ref{eq:v}), we can always assume that $d_{1}$ and $d_{2}$ are square-free. Thus, if $q|[d_{1},d_{2}]^{2},$ $q$ cannot be divided by the third power of a prime number. Therefore, we can write $q=q_{1}q_{2}^{2},$ where $(q_{1},q_{2})=1,$ $\mu(q_{1}q_{2})\ne 0.$ We note that
\begin{eqnarray}\label{eq:p0}
q|[d_{1},d_{2}]^{2}\Leftrightarrow q_{1}q_{2}|[d_{1},d_{2}],
\end{eqnarray} and
\begin{eqnarray}\label{eq:q1}
q_{1}q_{2}\geq q^{1/2}.
\end{eqnarray}
\begin{lemma}\label{le:1d} Set $U=N^{1/2-\delta/100}.$
For two co-prime integers $a$ and $q=q_{1}q_{2}^{2},$ where $(q_{1},q_{2})=1,$ $\mu(q_{1}q_{2})\ne 0,$ and $Y\leq N,$ there is
\begin{eqnarray*}
\sum\limits_{d_{1},d_{2}|P \atop q|[d_{1},d_{2}]^{2}}p_{d_{1}}p_{d_{2}}f_{d_{1},d_{2}}(Y,a/q)&=&\epsilon_{q}YT(q_{1}q_{2})+O\left(Y^{1/2}W^{1/2}z^{2}\right),
\end{eqnarray*}
where $T(q)$ is defined in (\ref{eq:tq}) and $\epsilon_{q}$ is as defined in Lemma \ref{le:1}.
\end{lemma}
{\it Proof of Lemma \ref{le:1d}:}
From (\ref{eq:tq}), Lemma \ref{le:1}, (\ref{eq:p0}) and the fact that $|p_{d}|\leq 1,$ we see
\begin{eqnarray}\label{eq:uhh}
&&\sum\limits_{d_{1},d_{2}|P \atop q|[d_{1},d_{2}]^{2}}p_{d_{1}}p_{d_{2}}f_{d_{1},d_{2}}(Y,a/q)=
\sum\limits_{d_{1},d_{2}|P \atop q_{1}q_{2}|[d_{1},d_{2}]}p_{d_{1}}p_{d_{2}}f_{d_{1},d_{2}}(Y,a/q)\nonumber\\
&=&\epsilon_{q}YT(q_{1}q_{2})+O\left(Y^{1/2}W^{1/2}\left(\sum\limits_{d\leq z}|p_{d}|\right)^{2}\right)
=\epsilon_{q}YT(q_{1}q_{2})+O\left(Y^{1/2}W^{1/2}z^{2}\right).\nonumber\\
\end{eqnarray}
\subsubsection{ Proof of Lemma \ref{le:pseudo1}} \label{se:ko}
In addition to the set of major arcs $\mathfrak{M }$ and minor arcs $\mathfrak{m} $ defined in (\ref{eq:arcs}), we defined a second set of major arcs $\mathfrak{M_{1}}$ and minor arcs $\mathfrak{m_{1}}$ as follows: We recall the definition of $R=\lfloor N^{1-\delta/2}\rfloor,$ set $U=N^{1/2-\delta/100},$ and define
\begin{eqnarray}\label{eq:arcs1}
\mathfrak{M_{1}}&=&\bigcup_{q\leq U}\bigcup_{a=1\atop (a,q)=1}^{q}\mathfrak{M}(a_{1},q_{1}),\nonumber\\ \mathfrak{M_{1}}(a,q)&=&\left\{\alpha\in [0,1]:\left|\alpha-\frac{a}{q}\right|\leq \frac{1}{qR}\right\},\nonumber\\
\mathfrak{m_{1}}&=&\mathbb{Z}/N\mathbb{Z}\setminus \mathfrak{M_{1}}.
\end{eqnarray} We prove Lemma \ref{le:pseudo1} separately for $r/N\in \mathfrak{M_{1}}$ and $r/N\in \mathfrak{m_{1}}$ in the next two paragraphs.
 \paragraph{The major arc case:} \label{p:2}
For $r/N\in \mathfrak{M_{1}},$ we write
\begin{eqnarray}\label{eq:popo}\frac{r}{N}=\frac{a}{q}+\beta,\quad q\leq U,\quad (a,q)=1,\quad |\beta|\leq 1/qR.\end{eqnarray} If $\beta=0,$ we see from Lemma \ref{le:1d}
\begin{eqnarray}\label{eq:b2u}
&&\frac{\phi(W)\log\,z}{W}\sum\limits_{d_{1},d_{2}|P\atop q| [d_{1},d_{2}]^{2}}p_{d_{1}}p_{d_{2}}f_{d_{1},d_{2}}(Y,a/q)\nonumber\\
&=&\frac{\phi(W)\log\,z}{W}\left(\epsilon_{q}YT(q_{1}q_{2})+E(Y)\right),
\end{eqnarray}
where
\begin{eqnarray*}
E(Y)&=&O\left( Y^{1/2}W^{1/2}z^{2}\right).
\end{eqnarray*}
Applying partial summation, we derive from (\ref{eq:defv1}) and (\ref{eq:b2u}),
\begin{eqnarray}\label{eq:M2aa}
\hat{v}_{1}\left(\frac{r}{N}\right)&=&\frac{\phi(W)\log\,z}{W}\int\limits_{0}^{N}e(\beta x)d\left(\sum\limits_{d_{1},d_{2}|P}p_{d_{1}}p_{d_{2}}f_{d_{1},d_{2}}(x,a/q)
\right)\nonumber\\
&=&\frac{\phi(W)\log\,z}{W}\epsilon_{q}T(q_{1}q_{2})\int\limits_{0}^{N}e(\beta x) dx
+\frac{\phi(W)\log\,z}{W}\int\limits_{0}^{N} e(\beta x) dE(x).\nonumber\\
\end{eqnarray}
Using (\ref{eq:popo}), we estimate the error term in (\ref{eq:M2aa}) as follows:
\begin{eqnarray}\label{eq:M3aa}
&&\left|\int\limits_{0}^{N} e(\beta x) dE(x)\right|\ll \left|E(N)\right|+
\left|\int\limits_{0}^{N}E(x)(2\pi i\beta) e(\beta )dx\right|\nonumber\\
&\ll & N^{1/2}W^{1/2}z^{2}(1+|\beta |N)\nonumber\\
&\ll & N^{1/2}W^{1/2}z^{2}\frac{N}{R}\nonumber\\
&\ll& N^{1-\delta/10}.
\end{eqnarray}
Combining (\ref{eq:M2aa}) -  (\ref{eq:M3aa}), we obtain
\begin{eqnarray}\label{eq:mi23aa}
\hat{v}_{1}\left(\frac{r}{N}\right)
&=&\frac{\phi(W)}{W}\log\,z \epsilon_{q}T(q_{1}q_{2})\int\limits_{0}^{N}e(\beta x) dx+O\left(N^{1-\delta/10}\right).\nonumber\\
\end{eqnarray} We now derive Lemma \ref{le:pseudo1} from (\ref{eq:mi23aa}): If $w<q\leq U,$ then Lemma \ref{le:pseudo1} follows from (\ref{eq:J}), Lemma \ref{le:a3}, and (\ref{eq:q1}). If $1<q\leq w,$ then $(q,W)>1,$ and thus $\epsilon_{q}=0.$
If $q=1$ and $\beta>0,$ then $\beta $ is an integer multiple of $1/N,$ and thus the integral in (\ref{eq:mi23aa}) equals zero. Finally, if $q=1$ and $\beta=0,$ then $\epsilon_{q}=1.$ Thus, using (\ref{eq:J}) and Lemma \ref{le:a3}, we obtain
\begin{eqnarray*}
\hat{v}_{1}\left(0\right)
&=&\frac{\phi(W)}{W}\log\,z\left(J^{-1}+O_{w}\left(w^{-1/2+\delta}\right)\right)N=(1+O_{w}(w^{-1/2+\delta})N.
\end{eqnarray*} This proves Lemma \ref{le:pseudo1} for the major arcs case sufficiently large $N$ and $z.$
 \paragraph{The minor arc case:} \label{p:3}
By (\ref{eq:arcs1}) and Dirichlet's theorem on rational approximation, we can write
\begin{eqnarray}\label{eq:u7y}\frac{r}{N}=\frac{a}{q}+\beta,\quad  U<q\leq R,\quad (a,q)=1,\quad |\beta|\leq q^{-2}.\end{eqnarray}
We argue as in the major arc case and derive (\ref{eq:M2aa}). For the major term in (\ref{eq:M2aa}), we argue in the same way as for the estimate of the major term in (\ref{eq:mi23aa})
in the case $w<q\leq U,$ and obtain an upper bound $\ll NU^{-1/4}.$
Estimating the error term as in (\ref{eq:M3aa}) and using (\ref{eq:u7y}), we see
\begin{eqnarray*}\label{eq:M3aabb}
&&\left|\int\limits_{0}^{N} e(\beta x) dE(x)\right|\ll \left|E(N)\right|+
\left|\int\limits_{0}^{N}E(x)(2\pi i\beta) e(\beta )dx\right|\nonumber\\
&\ll & N^{1/2}W^{1/2}z^{2}(1+|\beta |N)\nonumber\\
&\ll & N^{1/2}W^{1/2}z^{2}\frac{N}{U^{2}}\nonumber\\
&\ll& N^{1-\delta/10}.
\end{eqnarray*} This proves Lemma \ref{le:pseudo1} for the minor arcs case for sufficiently large $N$ and $z.$\\
\subsection{Proof of Lemma \ref{le:pseudo1a}}\label{se:a2}
We first assume $\frac{r}{N}=\frac{a}{q}$ and analyse the term of $f_{d_{1},d_{2}}(Y,a/q).$ For the analysis, we separately consider the two cases
1) $q_{d_{1},d_{2}}$ is $w$-smooth and/or $h\nmid 2$ and
2) $w$-rough and $h|2$ in  sections \ref{p:1} and \ref{p:2}. We then prove Lemma \ref{le:pseudo1a} in section \ref{se:ccc}.
\subsubsection{$q_{d_{1},d_{2}}$ is $w$-smooth and/or $h\nmid2$:}\label{p:1}
We will use the following Lemma:
\begin{lemma}\label{le:Sa} For $z\in [W]$ satisfying $z^{2}-b\equiv 0(mod\,W),$ we define
\begin{eqnarray}\label{eq:x1}
S_{q}(a,z)=\sum\limits_{r=1}^{q}e\left(\frac{a\left(Wr^{2}+2zr+\frac{z^{2}-b}{W}\right)}{q}\right).
\end{eqnarray}
a) If $(a,q)=1,$ $q>1,$ and $q$ is $w$-smooth, then
\begin{eqnarray*}
\sum\limits_{z\in [W]\atop z^{2}-b\equiv 0(mod\,W)}S_{q}(a,z)&=&0.
\end{eqnarray*}
b) If $(a,q)=1$ and $(q,W)\nmid 2,$
\begin{eqnarray*}
S_{q}(a,z)&=&0.
\end{eqnarray*}
c) For $(a,q)=1,$ there is,
\begin{eqnarray*}
|S_{q}(a,z)|&\leq &2\sqrt{q} .\end{eqnarray*}
\end{lemma} {\it Proof:} Part a) and c) are Lemma 5.3 and Lemma 5.2 in \cite{browning}, respectively. Part b) is stated in the proof of \cite[Lemma 5.3]{browning}.
\begin{lemma}\label{le:1d23}
 If $q_{d_{1},d_{2}}>1$ is $w$-smooth and/or $h\nmid 2,$ then for any two co-prime integers $a$ and $q$ there is:
\begin{eqnarray*}
f_{d_{1},d_{2}}(Y,a/q)\ll Y^{1/2}q_{d_{1},d_{2}}^{1/2}.
\end{eqnarray*}
\end{lemma}
{\it Proof of Lemma \ref{le:1d23}:}
We analyze the right-hand side of (\ref{eq:pseudo11111}).
As by assumption $q\nmid [d_{1}d_{2}]^{2},$ we have  $1< q_{d_{1},d_{2}}\leq q.$ We see
\begin{eqnarray}\label{eq:h66}
&&\sum\limits_{x\leq \sqrt{WY+b}\atop {[d_{1},d_{2}]|x\atop x^{2}-b\equiv 0(mod\,W)}}x
 e\left(\frac{a\left( x^{2}-b\right)}{qW}\right)=
e\left(\frac{-ab}{qW}\right)[d_{1},d_{2}] \sum\limits_{m\leq \sqrt{WY+b}/[d_{1},d_{2}]\atop m^{2}-b\overline{[d_{1},d_{2}]}^{2}\equiv 0(mod\,W)}m
 e\left(\frac{a_{d_{1},d_{2}}m^{2}}{q_{d_{1},d_{2}}W}\right),\nonumber\\
 \end{eqnarray}
where $[d_{1},d_{2}]\overline{[d_{1},d_{2}]}\equiv 1(mod\,W).$  Splitting the inner summation over $m$ into rest classes modulo $W,$ we find
\begin{eqnarray}\label{eq:h66a}
&&\sum\limits_{m\leq \sqrt{WY+b}/[d_{1},d_{2}]\atop m^{2}-b\overline{[d_{1},d_{2}]}^{2}\equiv 0(mod\,W)}m
 e\left(\frac{a_{d_{1},d_{2}}m^{2}}{q_{d_{1},d_{2}}W}\right)\nonumber\\
 &=&\sum\limits_{z\in[W]\atop z^{2}-b\overline{[d_{1},d_{2}]}^{2}\equiv 0(mod\,W)}\sum\limits_{z+yW\leq \sqrt{WY+b}/[d_{1},d_{2}]}(z+Wy)e\left(\frac{a_{d_{1},d_{2}}\left(z+yW\right)^{2}}{q_{d_{1},d_{2}}W}\right)\nonumber\\
 \end{eqnarray}
Splitting the summation over $y$ into rest classes modulo $q_{d_{1},d_{2}},$ we can write the inner sum in (\ref{eq:h66a}) as follows:
 \begin{eqnarray}\label{eq:h68}
&&\sum\limits_{z+yW\leq \sqrt{WY+b}/[d_{1},d_{2}]}(z+Wy)e\left(\frac{a_{d_{1},d_{2}}\left(z+yW\right)^{2}}{q_{d_{1},d_{2}}W}\right)\nonumber\\
 &=&\sum\limits_{r=1}^{q_{d_{1},d_{2}}}e\left(\frac{a_{d_{1},d_{2}}(z+Wr)^{2}}{q_{d_{1},d_{2}}W}\right)
\sum\limits_{z+Wr+Wq_{d_{1},d_{2}}s\leq \sqrt{WY+b}/[d_{1},d_{2}]}z+Wr+Wq_{d_{1},d_{2}}s.\nonumber\\
.\end{eqnarray}
We now evaluate the inner sum over $s$ in (\ref{eq:h68}):
\begin{eqnarray}\label{eq:h69}
&&\sum\limits_{z+Wr+Wq_{d_{1},d_{2}}s\leq \sqrt{WY+b}/[d_{1},d_{2}]}z+Wr+Wq_{d_{1},d_{2}}s\nonumber\\
&=&Wq_{d_{1},d_{2}}\sum\limits_{s\leq ((\sqrt{WY+b}/[d_{1},d_{2}])-z-Wr)/Wq_{d_{1},d_{2}}}s
+O\left((z+Wr)\sum\limits_{s\leq \sqrt{WY+b}/Wq_{d_{1},d_{2}}[d_{1},d_{2}]}1\right)\nonumber\\
&=&Wq_{d_{1},d_{2}}\sum\limits_{s\leq \sqrt{WY+b}/Wq_{d_{1},d_{2}}[d_{1},d_{2}]}s\nonumber\\
&+&O\left(Wq_{d_{1},d_{2}}\sum\limits_{((\sqrt{WY+b}/[d_{1},d_{2}])-z-Wr)/Wq_{d_{1},d_{2}}\leq s\leq \sqrt{WY+b}/Wq_{d_{1},d_{2}}[d_{1},d_{2}]}s\right)\nonumber\\
&+&O\left(Wq_{d_{1},d_{2}}\sum\limits_{s\leq \sqrt{WY+b}/Wq_{d_{1},d_{2}}[d_{1},d_{2}]}1\right)\nonumber\\
&=&Wq_{d_{1},d_{2}}\sum\limits_{s\leq \sqrt{WY+b}/Wq_{d_{1},d_{2}}[d_{1},d_{2}]}s
+O\left(\sqrt{WY}/[d_{1},d_{2}]\right).
\end{eqnarray}
Combining (\ref{eq:h66}) - (\ref{eq:h69}), see see
\begin{eqnarray}\label{eq:h71}
&&\sum\limits_{x\leq \sqrt{WY+b}\atop {[d_{1},d_{2}]|x\atop x^{2}-b\overline{[d_{1},d_{2}]}^{2}\equiv 0(mod\,W)}}x
 e\left(\frac{a\left( x^{2}-b\right)}{qW}\right)\nonumber\\
 &=&e\left(\frac{-ab}{qW}\right)[d_{1},d_{2}]\sum\limits_{z\in[W]\atop z^{2}-b=0(mod\,W)}\sum\limits_{r=1}^{q_{d_{1},d_{2}}}e\left(\frac{a_{d_{1},d_{2}}(z+Wr)^{2}}{q_{d_{1},d_{2}}W}\right)\nonumber\\
 &\times &\left(Wq_{d_{1},d_{2}}\sum\limits_{s\leq \sqrt{WY+b}/Wq_{d_{1},d_{2}}[d_{1},d_{2}]}s
+O\left(\sqrt{WY}/[d_{1},d_{2}]\right)\right)\nonumber\\
 &=&e\left(\frac{-ab}{qW}+\frac{a_{d_{1},d_{2}}b}{q_{d_{1},d_{2}}W}\right)[d_{1},d_{2}]\sum\limits_{z\in[W]\atop z^{2}-b\overline{[d_{1},d_{2}]}^{2}\equiv0(mod\,W)}S_{q_{d_{1},d_{2}}}(a_{d_{1},d_{2}},z)\nonumber\\
 &\times &\left(Wq_{d_{1},d_{2}}\sum\limits_{s\leq \sqrt{WY+b}/Wq_{d_{1},d_{2}}[d_{1},d_{2}]}s
+O\left(\sqrt{WY}/[d_{1},d_{2}]\right)\right)\nonumber\\
&=&e\left(\frac{-ab}{qW}+\frac{a_{d_{1},d_{2}}b}{q_{d_{1},d_{2}}W}\right)
\frac{WY+b}{2W[d_{1},d_{2}]q_{d_{1},d_{2}}}
\sum\limits_{z\in[W]\atop z^{2}-b\overline{[d_{1},d_{2}]}^{2}\equiv0(mod\,W)}S_{q_{d_{1},d_{2}}}(a_{d_{1},d_{2}},z)\nonumber\\
&+&O\left(\sqrt{WY}
\sigma(b)\max\limits_{z\in[W]\atop z^{2}-b\overline{[d_{1},d_{2}]}^{2}\equiv0(mod\,W)}\left|S_{q_{d_{1},d_{2}}}(a_{d_{1},d_{2}},z)\right|\right).
\end{eqnarray} Applying Lemma \ref{le:Sa} c), we estimate the $O$-term in (\ref{eq:h71}) as follows:
\begin{eqnarray}\label{eq:h72}
\sqrt{WY}
\sigma(b)\max\limits_{z\in[W]\atop z^{2}-b\overline{[d_{1},d_{2}]}^{2}\equiv0(mod\,W)}\left|S_{q_{d_{1},d_{2}}}(a_{d_{1},d_{2}},z)\right|\ll Y^{1/2}q^{1/2}_{d_{1},d_{2}}.\
\end{eqnarray}
Lemma \ref{le:1d23} now follows from (\ref{eq:pseudo11111}), (\ref{eq:h71}), Lemma \ref{le:Sa} a) and b), and (\ref{eq:h72}).
\begin{lemma}\label{le:1e}
For two co-prime integers $a$ and $q$ there is:
\begin{eqnarray*}
\sum\limits_{d_{1},d_{2}|P \atop{ q\nmid [d_{1},d_{2}]^{2},\,q_{W,d_{1},d_{2}}\leq Q \atop q_{d_{1},d_{2}}\,w-smooth\,and/or\, h\nmid 2}}p_{d_{1}}p_{d_{2}}f_{d_{1},d_{2}}(Y,a/q)\ll
z^{2}Y^{1/2}Q^{1/2}.
\end{eqnarray*}
\end{lemma}
{\it Proof of Lemma \ref{le:1e}:}  Using Lemma \ref{le:1d23} and noting that $q_{d_{1},d_{2}}\leq Wq_{W,d_{1},d_{2}},$ we see
\begin{eqnarray}\label{eq:h1h}
&&\sum\limits_{d_{1},d_{2}|P \atop{ q\nmid [d_{1},d_{2}]^{2},\,q_{W,d_{1},d_{2}}\leq Q \atop q_{d_{1},d_{2}}\,w-smooth\,and/or\, h\nmid 2}}p_{d_{1}}p_{d_{2}}f_{d_{1},d_{2}}(Y,a/q)\ll
\sum\limits_{d_{1},d_{2}|P \atop q_{W,d_{1},d_{2}}\leq Q}p_{d_{1}}p_{d_{2}}f_{d_{1},d_{2}}(Y,a/q)\nonumber\\
&\ll &\left(\sum\limits_{d\leq z}1\right)^{2}Y^{1/2}Q^{1/2}
\leq z^{2}Y^{1/2}Q^{1/2},\nonumber\\
 \end{eqnarray}
 qed.
 \subsubsection{$q_{d_{1},d_{2}}$ is $w$-rough and $h|2$:} \label{p:2}
We will first derive the auxiliary Lemmas \ref{le:hg} - \ref{le:uuu}. Subsequently we will prove the main Lemmas \ref{le:gaussapplied} and \ref{le:1d23y} of this paragraph.\\
\\
We will make use of the generalized Gauss sum $G(a,b,c)$ defined as follows:
\begin{eqnarray}\label{eq:x2}
G(a,b,c)=\sum\limits_{n=1}^{c}e\left(\frac{an^{2}+bn}{c}\right).
\end{eqnarray}
\begin{lemma}\label{le:hg} For co-prime integers $g$ and $h,$ $(h,2)=1,$ let $\left(\frac{g}{h}\right)$ denote the Jacobi symbol.
\begin{eqnarray*}
G(a,0,c)&=&\left(
                   \begin{array}{ll}
                     \left(\frac{a}{c}\right)\sqrt{c}, & c\equiv 1(mod\,4), \\
                                        \left(\frac{a}{c}\right)i\sqrt{c}, & c\equiv 3(mod\,4), \\
                                                            \end{array}
                 \right)
\end{eqnarray*}
\end{lemma}
{\it Proof of Lemma \ref{le:hg}:} See \cite[Theorem 1.5.2]{berndt}.

\begin{lemma}\label{le:hg6ggg}
For any integers $a,b$ and $c$ with $(c,2a)=1,$ we have
\begin{eqnarray*}
\qquad\qquad G(a,b,c)= e\left(\frac{-(\overline{2})^{2}\overline{a}b^{2}}{c}\right)\times\left(
                   \begin{array}{ll}
                     \left(\frac{a}{c}\right)\sqrt{c}, & c\equiv 1(mod\,4), \\
                                        \left(\frac{a}{c}\right)i\sqrt{c}, & c\equiv 3(mod\,4), \\
                                                            \end{array}
                 \right),
\end{eqnarray*}
where the integers $\overline{2}$ and $\overline{a}$ modulo $c$ are defined through the relations $\overline{2}2\equiv \overline{a}a\equiv 1(mod\,c).$\\
\end{lemma}
{\it Proof of Lemma \ref{le:hg6ggg}:}
\begin{eqnarray}\label{eq:kjkj}
&&G(a,b,c)=
e\left(\frac{-(\overline{2})^{2}\overline{a}b^{2}}{c}\right)\sum\limits_{n=1}^{c}e\left(\frac{a\left(n^{2}+\overline{a}bn+(\overline{2}\overline{a}b)^{2}\right)}{c}\right)\nonumber\\
&=&
e\left(\frac{-(\overline{2})^{2}\overline{a}b^{2}}{c}\right)\sum\limits_{n=1}^{c}e\left(\frac{a\left(n+\overline{2}\overline{a}b)^{2}\right)}{c}\right)=
e\left(\frac{-(\overline{2})^{2}\overline{a}b^{2}}{c}\right)\sum\limits_{n=1}^{c}e\left(\frac{an^{2}}{c}\right)\nonumber\\
&=& e\left(\frac{-(\overline{2})^{2}\overline{a}b^{2}}{c}\right)G(a,0,c).
\end{eqnarray} Applying Lemma \ref{le:hg} to (\ref{eq:kjkj}), we derive Lemma \ref{le:hg6ggg}.
\begin{lemma}\label{le:estimate}For fixed $k|q,\,k<q,$ there is
\begin{eqnarray*}
\left|\sum\limits_{d_{1},d_{2}\vert P\atop (q,[d_{1},d_{2}]^{2})=k}\frac{p_{d_{1}}p_{d_{2}}}{[d_{1},d_{2}]}\right|\ll d(q) J^{-1}k^{-1/2+\epsilon}.
\end{eqnarray*}
\end{lemma}
{\it Proof of Lemma \ref{le:estimate}:}
We write $\frac{q}{k}=\prod\limits_{i\leq M}p_{i}^{\alpha_{i}},$ where the $p_{i}$ are different prime numbers, $\alpha_{i}\in \mathbb{Z}^{+},$ and $M\leq d(p/k)$ is an integer depending on the prime decomposition of $\frac{q}{k}.$
Applying the inclusion-exclusion principle, we find
\begin{eqnarray}\label{eq:h84}
&&\sum\limits_{d_{1},d_{2}\vert P\atop (q,[d_{1},d_{2}]^{2})=k}\frac{p_{d_{1}}p_{d_{2}}}{[d_{1},d_{2}]}\nonumber\\&=&\left(
\sum\limits_{d_{1},d_{2}|P\atop k|[d_{1},d_{2}]^{2}}+\sum\limits_{1\leq j\leq M}(-1)^{j}\sum\limits_{\beta_{1}\leq 1}..\sum\limits_{\beta_{M}\leq 1\atop \beta_{1}+..\beta_{M}=j}
 \sum\limits_{d_{1},d_{2}|P\atop k\prod\limits_{i\leq M}p_{i}^{\beta_{i}}|[d_{1},d_{2}]^{2}}\right)\frac{p_{d_{1}}p_{d_{2}}}{[d_{1},d_{2}]}.\nonumber\\
\end{eqnarray}
From (\ref{eq:h84}), we see
\begin{eqnarray}\label{eq:h84a}
\left|\sum\limits_{d_{1},d_{2}\vert P\atop (q,[d_{1},d_{2}]^{2})=k}\frac{p_{d_{1}}p_{d_{2}}}{[d_{1},d_{2}]}\right|&\leq &
d(q)\max\limits_{j:k|j,\,j|q}\left|\sum\limits_{d_{1},d_{2}|P\atop j|[d_{1},d_{2}]^{2}}\frac{p_{d_{1}}p_{d_{2}}}{[d_{1},d_{2}]}\right|.
\end{eqnarray}
 Defining
$j=j_{1}j_{2}$ in the same way we have defined $q=q_{1}q_{2}$ in (\ref{eq:p0}) and applying Lemma \ref{le:a3}, we obtain from (\ref{eq:q1}) and (\ref{eq:h84a}):
\begin{eqnarray*}
\left|\sum\limits_{d_{1},d_{2}\vert P\atop (q,[d_{1},d_{2}]^{2})=k}\frac{p_{d_{1}}p_{d_{2}}}{[d_{1},d_{2}]}\right|&\ll &
d(q)J^{-1}\max\limits_{j:k|j,\,j|q}j^{-1/2+\epsilon}\leq d(q)J^{-1}k^{-1/2+\epsilon}.
\end{eqnarray*}
qed.
\begin{lemma}\label{le:uuu}
Let $\left(\frac{g}{h}\right)$ denote the Jacobi symbol. Let $a,b$ and $v$ be three strictly positive integers satisfying $(a,b)=1,$ $v|b,$  and $\left(\frac{b}{v},2\right)=1.$ Then, for any strictly positive integer $c$ with
$(c^{2},b)=v$, there is
\begin{eqnarray*}
\left(\frac{a\frac{c^{2}}{v}}{\frac{b}{v}}\right)=S_{a,b,v},
\end{eqnarray*}
where $S_{a,b,v}\in\{-1,1\}$ depends on $a,b,v,$ but is independent of $c.$
\end{lemma}
{\it Proof of Lemma \ref{le:uuu}:}
We write the prime decomposition of $v$ as $v=\prod\limits_{i\leq U}p_{i}^{\gamma_{i}},\gamma_{i}\geq 1.$
Further, we write the prime decomposition of $c^{2}$ as $c^{2}=d^{2}\prod\limits_{i\leq U}p_{i}^{2\beta_{i}},$ where $2\beta_{i}\geq \max(\gamma_{i},2),$\,  $(d,b)=1,\,p_{i}|b.$ Thus, $c^{2}/v=d^{2}\prod\limits_{i
\leq U\atop 2\beta_{i} >\gamma_{i}}p_{i}^{2\beta_{i}-\gamma_{i}}.$ We see,
\begin{eqnarray}\label{eq:v8}&&
\left(\frac{a\frac{c^{2}}{v}}{\frac{b}{v}}\right)=\left(\frac{a}{\frac{b}{v}}\right)\left(\frac{d}{\frac{b}{v}}\right)^{2}\prod\limits_{i\leq U}\left(\frac{p_{i}}{\frac{b}{v}}\right)^{2\beta_{i}-\gamma_{i}}\nonumber\\
&=&\left(\frac{a}{\frac{b}{v}}\right)\prod\limits_{i\leq U\atop 2\beta_{i}-\gamma_{i}\,is\,odd}\left(\frac{p_{i}}{\frac{b}{v}}\right)
=\left(\frac{a}{\frac{b}{v}}\right)\prod\limits_{i\leq U\atop \gamma_{i}\,is\,odd}\left(\frac{p_{i}}{\frac{b}{v}}\right).
\end{eqnarray} The last term in (\ref{eq:v8}) does not depend on $c$ which proves the Lemma.
\begin{lemma}\label{le:gaussapplied}
If $q_{d_{1},d_{2}}$ is $w$-rough and $h|2,$ then for any two co-prime integers $a$ and $q$ there is:
\begin{eqnarray*}
f_{d_{1},d_{2}}(Y,a/q)
&=&e\left(\frac{-ab}{qW}\right)e\left(\frac{g_{(q,[d_{1},d_{2}]^{2})}}{W_{1}}\right)
\frac{Y}{[d_{1},d_{2}]\sqrt{q_{d_{1},d_{2}}h}}
V_{q_{W,d_{1},d_{2}}}\left(\frac{W_{1}a_{d_{1},d_{2}}}{q_{W,d_{1},d_{2}}}\right)\nonumber\\
&
+&O\left(Y^{1/2}q_{d_{1},d_{2}}^{1/2}\right),
\end{eqnarray*}where $h,$ $W_{1},$ and $q_{W,d_{1},d_{2}}$ are as defined in (\ref{eq:w3w3w3}), $g_{(q,[d_{1},d_{2}]^{2})}$ is an integer modulo $W$ which for fixed $a,q,b,$ and $W$ depends on $(q,[d_{1},d_{2}]^{2})$ only, and
\begin{eqnarray}\label{eq:p9}
V_{q_{W,d_{1},d_{2}}}&=&\left(
                  \begin{array}{ll}
                    1, &  q_{W,d_{1},d_{2}}\equiv 1(mod\,4),\\
                     i, & q_{W,d_{1},d_{2}}\equiv 3(mod\,4).  \\
                   \end{array}
                 \right)
\end{eqnarray}
\end{lemma}
{\it Proof of Lemma \ref{le:gaussapplied}:}
By the definitions (\ref{eq:x1}) and (\ref{eq:x2}) and the assumption $h|2,$
\begin{eqnarray}\label{eq:x3}
S_{q_{d_{1},d_{2}}}(a_{d_{1},d_{2}},z) &=& e\left(\frac{(z^{2}-b)a_{d_{1},d_{2}}}{q_{d_{1},d_{2}}W}\right)G(W_{1}a_{d_{1},d_{2}},2za_{d_{1},d_{2}}/h,q_{W,d_{1},d_{2}}). \nonumber\\
\end{eqnarray}
Inserting (\ref{eq:x3}) in (\ref{eq:h71}) and using (\ref{eq:h72}), we find
\begin{eqnarray}\label{eq:h73}
&&\sum\limits_{x\leq \sqrt{WY+b}\atop {[d_{1},d_{2}]|x\atop x^{2}-b\equiv 0(mod\,W)}}x
 e\left(\frac{a\left( x^{2}-b\right)}{qW}\right)\nonumber\\
&=&e\left(\frac{-ab}{qW}\right)
\frac{WY+b}{2W[d_{1},d_{2}]q_{d_{1},d_{2}}}
\sum\limits_{z\in[W]\atop z^{2}-b\overline{[d_{1},d_{2}]}^{2}\equiv 0(mod\,W)}e\left(\frac{a_{d_{1},d_{2}}z^{2}}{q_{d_{1},d_{2}}W}\right)G(W_{1}a_{d_{1},d_{2}},2za_{d_{1},d_{2}}/h,q_{W,d_{1},d_{2}})\nonumber\\
&+&O\left(Y^{1/2}q_{d_{1},d_{2}}^{1/2}\right).
\end{eqnarray} We recall that we assume $h\in\{1,2\}$ which implies that $(q_{W,d_{1},d_{2}},2)=1.$ Applying Lemma \ref{le:hg6ggg} to (\ref{eq:h73}), we derive
\begin{eqnarray}\label{eq:h74}
&&\sum\limits_{x\leq \sqrt{WY+b}\atop {[d_{1},d_{2}]|x\atop x^{2}-b\equiv 0(mod\,W)}}x
 e\left(\frac{a\left( x^{2}-b\right)}{qW}\right)\nonumber\\
&=&e\left(\frac{-ab}{qW}\right)
\frac{(WY+b)\sqrt{q_{W,d_{1},d_{2}}}}{2W[d_{1},d_{2}]q_{d_{1},d_{2}}}
\sum\limits_{z\in[W]\atop z^{2}-b\overline{[d_{1},d_{2}]}^{2}\equiv 0(mod\,W)}e\left(\frac{a_{d_{1},d_{2}}z^{2}}{q_{d_{1},d_{2}}W}\right)
e\left(\frac{-(\overline{2})^{2}\overline{W_{1}a_{d_{1}d_{2}}}(2za_{d_{1},d_{2}}/h)^{2}}{q_{W,d_{1},d_{2}}}\right)\nonumber\\
&\times &V_{{q_{W,d_{1},d_{2}}}}\left(\frac{W_{1}a_{d_{1},d_{2}}}{q_{W,d_{1},d_{2}}}\right)
+O\left(Y^{1/2}q_{d_{1},d_{2}}^{1/2}\right)\nonumber\\
&=&e\left(\frac{-ab}{qW}\right)
\frac{WY+b}{2W[d_{1},d_{2}]\sqrt{q_{d_{1},d_{2}}h}}
\sum\limits_{z\in[W]\atop z^{2}-b\overline{[d_{1},d_{2}]}^{2}\equiv 0(mod\,W)}e\left(\frac{a_{d_{1},d_{2}}z^{2}}{q_{d_{1},d_{2}}W}\right)
e\left(\frac{-(\overline{2})^{2}\overline{W_{1}a_{d_{1}d_{2}}}(2za_{d_{1},d_{2}}/h)^{2}}{q_{W,d_{1},d_{2}}}\right)\nonumber\\
&\times &V_{{q_{W,d_{1},d_{2}}}}\left(\frac{W_{1}a_{d_{1},d_{2}}}{q_{W,d_{1},d_{2}}}\right)
+O\left(Y^{1/2}q_{d_{1},d_{2}}^{1/2}\right)\nonumber\\
&=&e\left(\frac{-ab}{qW}\right)
\frac{Y}{2[d_{1},d_{2}]\sqrt{q_{d_{1},d_{2}}h}}
\sum\limits_{z\in[W]\atop z^{2}-b\overline{[d_{1},d_{2}]}^{2}\equiv 0(mod\,W)}e\left(\frac{a_{d_{1},d_{2}}z^{2}}{q_{d_{1},d_{2}}W}\right)
e\left(\frac{-(\overline{2})^{2}\overline{W_{1}a_{d_{1}d_{2}}}(2za_{d_{1},d_{2}}/h)^{2}}{q_{W,d_{1},d_{2}}}\right)\nonumber\\
&\times &V_{{q_{W,d_{1},d_{2}}}}\left(\frac{W_{1}a_{d_{1},d_{2}}}{q_{W,d_{1},d_{2}}}\right)
+O\left(Y^{1/2}q_{d_{1},d_{2}}^{1/2}\right).
\end{eqnarray}
We recall that in (\ref{eq:h74}) for any integer $c$ prime to $q_{W,d_{1},d_{2}},$ $\overline{c}$ is defined via the relation $c\overline{c}\equiv 1(mod\,q_{W,d_{1},d_{2}}).$
Thus, in particular
\begin{eqnarray}\label{eq:t2g}a_{d_{1},d_{2}}\overline{a_{d_{1},d_{2}}}\equiv 1(mod\,q_{W,d_{1},d_{2}}).\end{eqnarray}
Further, we set
\begin{eqnarray}\label{eq:t2}W_{1}\overline{W_{1}}=1+s_{(q,[d_{1},d_{2}]^{2})}q_{W,d_{1},d_{2}}.\end{eqnarray}
In view of the definitions in (\ref{eq:w3w3w3}), we see that - for fixed $q$ and $W$ - $s_{(q,[d_{1},d_{2}]^{2})}$ depends on $(q,[d_{1},d_{2}]^{2})$ only.
Similarly, we notice that for fixed values of $q$ and $W,$ the value of the integer $\overline{2}$ only depends on $(q,[d_{1},d_{2}]^{2}) ,$ i.e.,
\begin{eqnarray}\label{eq:t3}\overline{2}:=u_{(q,[d_{1},d_{2}]^{2})} .\end{eqnarray}
Finally, we consider a fixed $z$ satisfying the congruence condition in (\ref{eq:h74}), i.e., $z^{2}[d_{1},d_{2}]^{2}\equiv b(mod\,W).$ This implies that
\begin{eqnarray}\label{eq:t1}
a_{d_{1},d_{2}}z^{2}=\frac{a[d_{1},d_{2}]^{2}}{(q,[d_{1},d_{2}]^{2})}z^{2}\equiv ab\overline{(q,[d_{1},d_{2}]^{2})}(mod\,W),
\end{eqnarray}
where $(q,[d_{1},d_{2}]^{2})\overline{(q,[d_{1},d_{2}]^{2})}\equiv 1(mod\,W).$
Subsequently applying (\ref{eq:t2g}) - (\ref{eq:t1}), we can calculate the product of the exponential terms in
(\ref{eq:h74}) as follows:
\begin{eqnarray}\label{eq:h75}
&&e\left(\frac{a_{d_{1},d_{2}}z^{2}}{q_{d_{1},d_{2}}W}\right)
e\left(\frac{-(\overline{2})^{2}\overline{W_{1}a_{d_{1}d_{2}}}(2za_{d_{1},d_{2}}/h)^{2}}{q_{W,d_{1},d_{2}}}\right)\nonumber\\
&=&e\left(\frac{a_{d_{1},d_{2}}z^{2}}{q_{d_{1},d_{2}}W}\right)
e\left(\frac{-(\overline{2})^{2}W_{1}\overline{W_{1}}a_{d_{1}d_{2}}(2z/h)^{2}}{q_{W,d_{1},d_{2}}W_{1}}\right)
 \nonumber\\
  &=&e\left(\frac{a_{d_{1},d_{2}}z^{2}}{q_{d_{1},d_{2}}W}\right)
e\left(\frac{-a_{d_{1}d_{2}}z^{2}(\overline{2})^{2}(2/h)^{2}}{q_{W,d_{1},d_{2}}W_{1}}\right)e\left(\frac{-s_{(q,[d_{1},d_{2}]^{2})}a_{d_{1},d_{2}}z^{2}(\overline{2})^{2}(2/h)^{2}}{W_{1}}\right)
 \nonumber\\
 &=& e\left(\frac{a_{d_{1},d_{2}}z^{2}}{q_{d_{1},d_{2}}W}\right)e\left(\frac{-a_{d_{1},d_{2}}z^{2}}{q_{W,d_{1},d_{2}}W_{1}h^{2}}\right)
e\left(\frac{-s_{(q,[d_{1},d_{2}]^{2})}ab\overline{(q,[d_{1},d_{2}]^{2})}(u_{(q,[d_{1},d_{2}]^{2})})^{2}(2/h)^{2}}{W_{1}}\right)\nonumber\\
&=& e\left(\frac{a_{d_{1},d_{2}}z^{2}}{q_{d_{1},d_{2}}W}\right)e\left(\frac{-a_{d_{1},d_{2}}z^{2}}{q_{d_{1},d_{2}}W}\right)
e\left(\frac{-s_{(q,[d_{1},d_{2}]^{2})}ab\overline{(q,[d_{1},d_{2}]^{2})}(u_{(q,[d_{1},d_{2}]^{2})})^{2}(2/h)^{2}}{W_{1}}\right)\nonumber\\
&:=&e\left(\frac{g_{(q,[d_{1},d_{2}]^{2})}}{W_{1}}\right),
\end{eqnarray} where - for fixed $a,q,b,$ and $W$ - $g_{(q,[d_{1},d_{2}]^{2})}$ is an integer modulo $W$ that depends on $(q,[d_{1},d_{2}]^{2})$ only.
Inserting (\ref{eq:h75}) into (\ref{eq:h74}), we obtain
\begin{eqnarray}\label{eq:h77}
&&\sum\limits_{x\leq \sqrt{WY+b}\atop {[d_{1},d_{2}]|x\atop x^{2}-b\equiv 0(mod\,W)}}x
 e\left(\frac{a\left( x^{2}-b\right)}{qW}\right)
=\frac{\sigma(b)}{2}e\left(\frac{-ab}{qW}\right)e\left(\frac{g_{(q,[d_{1},d_{2}]^{2})}}{W_{1}}\right)
\frac{Y}{[d_{1},d_{2}]\sqrt{q_{d_{1},d_{2}}h}}\nonumber\\
&\times & V_{{q_{W,d_{1},d_{2}}}}\left(\frac{W_{1}a_{d_{1},d_{2}}}{q_{W,d_{1},d_{2}}}\right)
+O\left(Y^{1/2}q_{d_{1},d_{2}}^{1/2}\right)).
\end{eqnarray}
Now the Lemma follows from (\ref{eq:pseudo11111}) and (\ref{eq:h77}).
\begin{lemma}\label{le:1d23y}
For two co-prime integers $a$ and $q$ there is:
\begin{eqnarray*}
&&\sum\limits_{d_{1},d_{2}|P \atop{ q\not|[d_{1},d_{2}]^{2},\,q_{W,d_{1},d_{2}}\leq Q  \atop q_{d_{1},d_{2}}\,w-rough\,and\, h|2}}p_{d_{1}}p_{d_{2}}f_{d_{1},d_{2}}(Y,a/q)\nonumber\\
&=& \frac{Y}{\sqrt{h}}\sum\limits_{k< \frac{q}{w}.\,q/kh\leq Q\atop {k|q\atop \frac{q}{k} w-rough\,and\, (W,\frac{q}{k})|2}}\frac{t_{k}}{\sqrt{q/k}}\sum\limits_{d_{1},d_{2}|P \atop (q,[d_{1},d_{2}]^{2})=k }
\frac{p_{d_{1}}p_{d_{2}}}{[d_{1},d_{2}]}
+O\left(z^{2}Y^{1/2}Q^{1/2}\right),
\end{eqnarray*}
where $t_{k}$ is a complex number that for fixed $a,q,b,$ and $W$ depends on $k$ only and $|t_{k}|=1.$
\end{lemma}
{\it Proof of Lemma \ref{le:1d23y}:} Using Lemma \ref{le:gaussapplied}, we see
\begin{eqnarray}\label{eq:w4545bcb}
&&\sum\limits_{d_{1},d_{2}|P \atop{ q\not|[d_{1},d_{2}]^{2},\,q_{W,d_{1},d_{2}}\leq Q  \atop q_{d_{1},d_{2}}\,w-rough\,and\, h|2}}p_{d_{1}}p_{d_{2}}f_{d_{1},d_{2}}(Y,a/q)\nonumber\\
&=& e\left(\frac{-ab}{qW}\right)\frac{Y}{\sqrt{h}}\sum\limits_{d_{1},d_{2}|P \atop{ q\not|[d_{1},d_{2}]^{2},\,q_{W,d_{1},d_{2}}\leq Q \atop q_{d_{1},d_{2}}\,w-rough\,and\, (W,q_{d_{1}d_{2}})|2}}
\frac{p_{d_{1}}p_{d_{2}}}{[d_{1},d_{2}]\sqrt{q_{d_{1},d_{2}}}}e\left(\frac{g_{(q,[d_{1},d_{2}]^{2})}}{W_{1}}\right)
V_{{q_{W,d_{1},d_{2}}}}\left(\frac{W_{1}a_{d_{1},d_{2}}}{q_{W,d_{1},d_{2}}}\right)\nonumber\\
&+&O\left(z^{2}Y^{1/2}Q^{1/2}\right),
\end{eqnarray} where we have estimated the $O$- term arguing similarly as in (\ref{eq:h1h}). If $q_{d_{1},d_{2}}$ is $w$-rough, then
$\label{eq:yh}q_{d_{1},d_{2}}>w,$
which implies that $(q,[d_{1},d_{2}]^{2})<q/w.$ Thus, we can rewrite the right-hand side of (\ref{eq:w4545bcb}) as follows:
\begin{eqnarray}\label{eq:w4}
&&\sum\limits_{d_{1},d_{2}|P \atop{ q\not|[d_{1},d_{2}]^{2},\,q_{W,d_{1},d_{2}}\leq Q  \atop q_{d_{1},d_{2}}\,w-rough\,and\, h|2}}p_{d_{1}}p_{d_{2}}f_{d_{1},d_{2}}(Y,a/q)\nonumber\\
&=& e\left(\frac{-ab}{qW}\right)\frac{Y}{\sqrt{h}}\sum\limits_{k< \frac{q}{w},\,q/kh\leq Q\atop {k|q\atop \frac{q}{k} w-rough\,and\, (W,\frac{q}{k})|2}}\frac{e(g_{k}/W_{1})V_{\frac{q}{hk}}}{\sqrt{q/k}}\sum\limits_{d_{1},d_{2}|P \atop (q,[d_{1},d_{2}]^{2})=k }
\frac{p_{d_{1}}p_{d_{2}}}{[d_{1},d_{2}]}
\left(\frac{W_{1}a[d_{1},d_{2}]^{2}/k}{q/hk}\right)\nonumber\\
&+&O\left(z^{2}Y^{1/2}Q^{1/2}\right).
\end{eqnarray}
Applying Lemma \ref{le:uuu} with $a=W_{1}a,$ $b=q/h,$ $v=k,$ and $c=[d_{1},d_{2}],$ we see that
  $\left(\frac{W_{1}a[d_{1},d_{2}]^{2}/k}{q/hk}\right)=S_{W_{1}a,q/h,k}.$ Thus, for fixed $a,q$ and $W,$ $ \left(\frac{W_{1}a[d_{1},d_{2}]^{2}/k}{q/hk}\right)$ depends on $k$ only.
We now define the complex number $t_{k}$ as\\ $t_{k}:=e\left(\frac{-ab}{qW}\right)e(g_{k}/W_{1})V_{\frac{q}{hk}}\left(\frac{W_{1}a[d_{1},d_{2}]^{2}/k}{q/hk}\right).$ By the foregoing discussion and the definition of $g_{k}$ and $V_{a,W,\frac{q}{hk}},$ we see that for fixed $a,q,b,$ and $W,$ $t_{k},$ depends on $k$ only and $|t_{k}|=1.$ Thus, we can rewrite (\ref{eq:w4}) as
\begin{eqnarray*}
&&\sum\limits_{d_{1},d_{2}|P \atop{ q\not|[d_{1},d_{2}]^{2} ,\,q_{W,d_{1},d_{2}}\leq Q \atop q_{d_{1},d_{2}}\,w-rough\,and\, h|2}}p_{d_{1}}p_{d_{2}}f_{d_{1},d_{2}}(Y,a/q)\nonumber\\
&=& \frac{Y}{\sqrt{h}}\sum\limits_{k< \frac{q}{w},\,q/kh\leq Q\atop {k|q\atop \frac{q}{k} w-rough\,and\, (W,\frac{q}{k})|2}}\frac{t_{k}}{\sqrt{q/k}}\sum\limits_{d_{1},d_{2}|P \atop (q,[d_{1},d_{2}]^{2})=k }
\frac{p_{d_{1}}p_{d_{2}}}{[d_{1},d_{2}]}
+\left(z^{2}Y^{1/2}Q^{1/2}\right),
\end{eqnarray*} qed.
\subsubsection{Proof of Lemma \ref{le:pseudo1a}}\label{se:ccc}
By (\ref{eq:arcs}) and (\ref{eq:defv1}), for fixed $r,$ we only need to consider those pairs $d:=(d_{1},d_{2})$ for which
\begin{eqnarray*}\label{eq:m8888}\left|\frac{rW[d_{1},d_{2}]^{2}}{N}-\frac{a_{1d}}{q_{1d}}\right|\leq \frac{1}{q_{1d}R},\quad\, (a_{1d},q_{1d})=1,\quad q_{1d}\leq Q,\end{eqnarray*} which implies
\begin{eqnarray}\label{eq:mumumu}&&\frac{r}{N}=\frac{a_{d}}{q_{d}}+\beta_{d},\quad  q_{d}\leq q_{1d}W[d_{1},d_{2}]^{2}\leq QWz^{2}\leq N^{1/2-\delta/100}, \quad |\beta_{d}|\leq 1/R.\nonumber\\
\end{eqnarray}
We now show that for different pairs $d:=(d_{1},d_{2})$ and $d^{*}:=(d^{*}_{1},d^{*}_{2})$ which satisfy (\ref{eq:mumumu}), there is
$a_{d}=a^{*}_{d}$ and $q_{d}=q^{*}_{d}$ which implies that $\beta_{d}=\beta_{d^{*}}.$ This follows from (\ref{eq:mumumu}) and the relation
\begin{eqnarray*}
\left|\frac{a_{d}}{q_{d}}-\frac{a^{*}_{d}}{q^{*}_{d}}\right|\geq \frac{1}{q_{d}q^{*}_{d}}>\frac{1}{q_{d}R}+\frac{1}{q_{d^{*}}R},\end{eqnarray*} which holds because of
$q_{d},q_{d^{*}}\leq N^{1/2-\delta/100}.$
Thus, we rewrite (\ref{eq:mumumu}) as
\begin{eqnarray}\label{eq:mumumu1a}&&\frac{r}{N}=\frac{a}{q}+\beta,\quad (a,q)=1,\quad q\leq N^{1/2-\delta/100},\quad  |\beta|\leq 1/R.\end{eqnarray}
In view of (\ref{eq:mumumu1a}), we first consider the case $\beta=0,$ i.e., $\frac{r}{N}=\frac{a}{q}.$
By (\ref{eq:arcs}) and (\ref{eq:defv1}), for fixed $r,$ we see that $rW[d_{1},d_{2}]^{2}/N\in \mathfrak{M}$ if and only if
\begin{eqnarray}\label{eq:m8888}\frac{rW[d_{1},d_{2}]^{2}}{N}=\frac{\frac{aW[d_{1},d_{2}]^{2}}{(q,W[d_{1},d_{2}]^{2})}}{q_{W,d_{1},d_{2}}},\quad q_{W,d_{1},d_{2}}\leq Q.
\end{eqnarray}
From (\ref{eq:defv1}), Lemma \ref{le:1e}, Lemma \ref{le:1d23y}, and (\ref{eq:m8888}), we see
\begin{eqnarray}\label{eq:b2} \hat{v}_{2}\left(\frac{a}{q}\right)&=&
\frac{\phi(W)\log\,z}{W}\sum\limits_{d_{1},d_{2}|P\atop{q\nmid [d_{1},d_{2}]^{2}\atop q_{W,d_{1},d_{2}}\leq Q}}p_{d_{1}}p_{d_{2}}f_{d_{1},d_{2}}(Y,a/q)\nonumber\\
&=&\frac{\phi(W)\log\,z}{W}\left(YM(q,a)+E_{1}(Y)\right),
\end{eqnarray}
where
\begin{eqnarray*}
M(a,q)&=&\frac{1}{\sqrt{h}}\sum\limits_{k< \frac{q}{w},\,q/kh\leq Q\atop {k|q\atop \frac{q}{k} w-rough\,and\, (W,\frac{q}{k})|2}}\frac{t_{k}}{\sqrt{q/k}}\sum\limits_{d_{1},d_{2}|P \atop (q,[d_{1},d_{2}]^{2})=k }
\frac{p_{d_{1}}p_{d_{2}}}{[d_{1},d_{2}]},\\
E_{1}(Y)&=&O\left( z^{2}Y^{1/2}Q^{1/2}\right).
\end{eqnarray*}
We notice from the definition of $M(a,q)$ that
\begin{eqnarray}\label{eq:p999}
M(a,q)\ne 0\Rightarrow q>w.
\end{eqnarray}
We now consider the general case $\frac{r}{N}=\frac{a}{q}+\beta,$  Applying partial summation, we derive from (\ref{eq:defv1}), (\ref{eq:mumumu1a}), and (\ref{eq:b2}),
\begin{eqnarray}\label{eq:M2u}
\hat{v}_{2}\left(\frac{r}{N}\right)
&=&\frac{\phi(W)\log\,z}{W}M(q,a)\int\limits_{0}^{N}e(\beta x) dx\nonumber\\
&+&\frac{\phi(W)\log\,z}{W}\int\limits_{0}^{N} e(\beta x) dE_{1}(x).
\end{eqnarray}
We first estimate the main term in (\ref{eq:M2u}). Applying Lemma \ref{le:estimate} and using (\ref{eq:J}), we find:
\begin{eqnarray}\label{eq:M31uaaa}
&&\frac{\phi(W)\log\,z}{W}M(q,a)\ll  \log z\,d(q)\max\limits_{k| q,\,k<q/w}(q/k)^{-1/2} \left|\sum\limits_{d_{1},d_{2}|P \atop (q,[d_{1},d_{2}]^{2})=k }
\frac{p_{d_{1}}p_{d_{2}}}{[d_{1},d_{2}]}\right|\nonumber\\
& \ll& \log z\, J^{-1}\max\limits_{k| q,\,k<q/w}(q/k)^{-1/2+\epsilon}k^{-1/2+\epsilon} \ll  q^{-1/2+2\epsilon}.
\end{eqnarray} We recall the well-known estimate
\begin{eqnarray}\label{eq:M31uab}
\int\limits_{0}^{N}e(\beta x) dx\ll \min\left(N,\frac{1}{||\beta ||}\right).
\end{eqnarray}
Combining (\ref{eq:p999}), (\ref{eq:M31uaaa}), and (\ref{eq:M31uab}), we find
\begin{eqnarray}\label{eq:M31u}
&&\frac{\phi(W)\log\,z}{W}M(q,a)\int\limits_{0}^{N}e(\beta x) dx\ll Nw^{-1/2+2\epsilon}.
\end{eqnarray}
Using (\ref{eq:mumumu1a}), we estimate the error term integral in (\ref{eq:M2u}) as follows:
\begin{eqnarray}\label{eq:M3u}
&&\left|\int\limits_{0}^{N} e(\beta x) dE_{1}(x)\right|\ll \left|E_{1}(N)\right|+
\left|\int\limits_{0}^{N}E_{1}(x)(2\pi i\beta) e(\beta )dx\right|\nonumber\\
&\ll & \left(z^{2}N^{1/2}Q^{1/2}\right)(1+|\beta |N)\nonumber\\
&\ll & \left(z^{2}N^{1/2}Q^{1/2}\right)\frac{N}{R}\nonumber\\
&\ll& N^{1-\delta/10}.
\end{eqnarray}
Lemma \ref{le:pseudo1a} now follows from (\ref{eq:M2u}), (\ref{eq:M31u}), and (\ref{eq:M3u}).
\subsection{Proof of Lemma \ref{le:pseudo2}}
Using (\ref{eq:pseudo11111}) with $\frac{r}{N}$ instead of $\frac{a}{q},$ we see
\begin{eqnarray}\label{eq:mi1}
&&\left|f_{d_{1},d_{2}}(N,r/N)\right|= \left|\frac{2}{\sigma(b)}\sum\limits_{x\leq \sqrt{WN+b}\atop {[d_{1},d_{2}]|x\atop x^{2}-b\equiv 0(mod\,W)}}x
 e\left(\frac{r\left( x^{2}-b\right)}{NW}\right)\right|\nonumber\\
&\ll & \sum\limits_{j=1}^{\sigma(b)}\left|
\sum\limits_{x\leq \sqrt{WN+b}\atop {[d_{1},d_{2}]|x\atop x\equiv  h_{j}(mod\,W)}}x
 e\left(\frac{r x^{2}}{NW} \right)\right|.\nonumber\\
&=& \sum\limits_{j=1}^{\sigma(b)}\left| \sum\limits_{x\leq \sqrt{WN+b} \atop x\equiv  g_{j}(mod\,W[d_{1},d_{2}])} xe\left(\frac{r x^{2}}{NW}\right)\right|,
 \end{eqnarray}
 where $g_{j}=g_{j}(h_{j},W[d_{1},d_{2}]).$ Here we have used the fact that due to $(W,P)=1$ and $d_{1},d_{2}|P,$ there is $(W,[d_{1},d_{2}])=1.$
 Using partial summation, we estimate the inner sum over $x$ in (\ref{eq:mi1}) for $Y\leq N$ as follows:
 \begin{eqnarray}\label{eq:h2}
\left|
 \sum\limits_{x\leq \sqrt{WN+b} \atop x\equiv  g_{j}(mod\,W[d_{1},d_{2}])} xe\left(\frac{r x^{2}}{NW}\right)\right|
&\ll & \sqrt{WN}\max_{K\leq \sqrt{WN+b}} \left|U(K,r/N)\right|,\nonumber\\
\end{eqnarray}
where
\begin{eqnarray}\label{eq:h211}U(K,r/N)&:=&\sum\limits_{x\leq K \atop x\equiv  g_{j}(mod\,W[d_{1},d_{2}])} e\left(\frac{r x^{2}}{NW}\right)\nonumber\\
&=& \sum\limits_{s\leq (K-g_{j})/W[d_{1},d_{2}]}e\left(\frac{r W[d_{1},d_{2}]^{2}s^{2}}{N} +\frac{2r g_{g}[d_{1},d_{2}] s}{N} + \frac{r g_{j}^{2}}{NW} \right)\nonumber\\.
\end{eqnarray}
To estimate the right-hand side of (\ref{eq:h211}), we wil make
use of Weyl's Lemma \cite[Lemma 2.4]{vaughan}:
\begin{lemma}\label{le:weyl} Let $\alpha, $ $\alpha_{1}$ and $\alpha_{2}$ be real numbers. If
$\left|\alpha-\frac{a}{q}\right|\leq q^{-2},$ for two integers $a$ and $q$ with $(a,q)=1,$ then
\begin{eqnarray*}
\sum\limits_{n\leq N}e\left(\alpha n^{2}+\alpha_{1}n+\alpha_{2}\right)\ll N^{1+\epsilon}\left(q^{-1}+N^{-1}+qN^{-2}\right)^{1/2}.
\end{eqnarray*}
\end{lemma}
By the definition of the minor arcs (\ref{eq:arcs}) and Dirichlet's theorem on rational approximation, we know that there exists integers $a$ and $q$ with $(a,q)=1,$ and $Q<q\leq R$ such that
$\left|\frac{rW[d_{2},d_{2}]^{2}}{N}-\frac{a}{q}\right|\leq 1/q^{2}.$
Thus, applying Lemma \ref{le:weyl} to (\ref{eq:h211}) with $\alpha =rW[d_{2},d_{2}]^{2}/N,$ we find
\begin{eqnarray}\label{eq:mi2b}
 U(K,r/N)&\ll &\frac{K^{1+\epsilon}}{W^{1+\epsilon}[d_{1},d_{2}]^{1+\epsilon}}\left(Q^{-1}+K^{-1}+RK^{-2}\right)^{1/2}.
\end{eqnarray} The right-hand side of (\ref{eq:mi2b}) is an increasing function in $K.$ Therefore,
\begin{eqnarray}\label{eq:mi2c}
&&\max_{K\leq \sqrt{WN+b}}\left| U\left(K,r/N\right)\right|\nonumber\\
&\ll &\frac{N^{1/2+\epsilon}}{[d_{1},d_{2}]^{1+\epsilon}}
\left(Q^{-1}+N^{-1/2}+RN^{-1}\right)^{1/2}.\nonumber\\
&\ll & \frac{N^{1/2+\epsilon}}{[d_{1},d_{2}]^{1+\epsilon}}
\left( N^{-2\delta/5}+N^{-1/2}+N^{-\delta/2}\right)^{1/2}
\nonumber\\
&\ll &
\frac{N^{\frac{1}{2}-\frac{\delta}{100}}}{[d_{1},d_{2}]}.\end{eqnarray}
From (\ref{eq:h2}) and (\ref{eq:mi2c}), we see
\begin{eqnarray}\label{eq:mi2d}
\left| \sum\limits_{x\leq \sqrt{WN+b} \atop x\equiv  g_{j}(mod\,W[d_{1},d_{2}])} xe\left(\alpha x^{2}/W\right)\right|&\ll & \frac{N^{1-\frac{\delta}{100}}}{[d_{1},d_{2}]}.
\end{eqnarray}
From (\ref{eq:defv1}), (\ref{eq:mi1}), (\ref{eq:mi2d}) and using $|\rho_{d}|\leq 1,$ we see
\begin{eqnarray}\label{eq:mi2eccd}
\hat{v}_{3}\left(\frac{r}{N}\right) &\ll& N^{1-\frac{\delta}{150}}\sum\limits_{d_{1},d_{2}|P}\frac{1}{[d_{1},d_{2}]}\ll N^{1-\frac{\delta}{150}}\sum\limits_{d_{1},d_{2}\leq z}\frac{1}{[d_{1},d_{2}]}\nonumber\\
&\ll & N^{1-\frac{\delta}{150}}\sum\limits_{d\leq z^{2}\atop (d,W)=1}\left(\sum\limits_{d_{1},d_{2}\leq z\atop [d_{1},d_{2}]=d}1\right)\frac{1}{d}.
\end{eqnarray}
For any fixed squarefree $d\leq z^{2},$ there are at most $3^{\varpi(d)}\leq d^{\delta/900}\leq  N^{\delta/900}$ pairs $[d_{1},d_{2}]=d.$ Thus, we see from (\ref{eq:mi2eccd})
\begin{eqnarray*}
\hat{v}_{3}\left(\frac{r}{N}\right)&\ll &N^{1-\delta/300}\sum\limits_{d\leq z^{2}}d^{-1}\ll \ll N^{1-\delta/400},
\end{eqnarray*}
which proves Lemma \ref{le:pseudo2}.

\section{Restriction estimate}\setcounter{equation}{0}\setcounter{theorem}{0}\setcounter{lemma}{0}
In this section, we set $a(n)=a_{i}(n),$ $b=b_{i},$ $z=z_{i},$ $v(n)=v_{i}(n),$ and $N=N_{i}$ The main purpose of this section is to show the following Lemma \ref{le:r1}. Our proof follows the argument in \cite[Section 6]{browning} with some minor modifications.
\begin{lemma}\label{le:r1}
        For any real number $p>4$ there exists an absolute constant $C_{p}$ such that
        \begin{eqnarray*}
\int\limits_{\mathbb{T}}\left|\hat{a}(\theta)\right|^{p}d\theta \leq C_{p}N^{p-1}.
           \end{eqnarray*}
\end{lemma}
For the proof of Lemma \ref{le:r1}, we will make use of the following Lemma:
\begin{lemma}\label{le:r2}
        There exists an absolute constant $C$ such that
        \begin{eqnarray*}
\int\limits_{\mathbb{T}}\left|\hat{a}(\theta)\right|^{4}d\theta \leq N^{3+C/\log\log N}.
           \end{eqnarray*}
\end{lemma}
{\it Proof of Lemma \ref{le:r2}:}
We note that
\begin{eqnarray}\label{eq:rw}
&&\int\limits_{\mathbb{T}}\left|\hat{a}(\theta)\right|^{4}d\theta=\left|\left|\sum\limits_{n,m\in [1,N]}a(m)\bar{a}(n)e(\theta(m-n))\right|\right|_{2}^{2}\nonumber\\
&\leq & \sum\limits_{|k|\leq N}
\left|\sum\limits_{n,m\in [1,N]\atop n-m=k}a(m)\bar{a}(n)\right|^{2}.\end{eqnarray}
The contribution to the sum over $k$ for $k=0$ is
\begin{eqnarray}\label{eq:rw1}
\leq \left|\sum\limits_{n\in [1,N]}|a(n)|^{2}\right|^{2}\ll N^{3}(\log\,N)^{2}=N^{3+O(1/\log\log\,N)}.
\end{eqnarray}
If $k\ne 0,$ we see from the definition of $a(n),$
\begin{eqnarray}\label{eq:MMMaaa}
&&\sum\limits_{n,m\in [1,N]\atop n-m=k}a(m)\bar{a}(n)=\left(\frac{2\phi(W)\log\,z}{\sigma(b)W}\right)^{2}\sum\limits_{z\leq x,y\leq \sqrt{WN+b},\,x,y\,prime\atop {x^{2}\equiv y^{2}\equiv b(mod\,W)\atop  x^{2}-y^{2}=Wk}}xy.
\end{eqnarray} Using the Cauchy-Schwarz inequality, we estimate the sum on the right hand side in (\ref{eq:MMMaaa}) as follows:
\begin{eqnarray}\label{eq:MMMaaaaa}
\sum\limits_{z\leq x,y\leq \sqrt{WN+b},\,x,y\,prime\atop {x^{2}\equiv y^{2}\equiv b(mod\,W)\atop  x^{2}-y^{2}=Wk}}xy\ll d^{1/2}(k)
\left(\sum\limits_{x,y\leq \sqrt{WN+b},\,x,y\,prime\atop {x^{2}\equiv y^{2}\equiv b(mod\,W)\atop  x^{2}-y^{2}=Wk}}x^{2}y^{2}\right)^{1/2}.
\end{eqnarray}
Combining (\ref{eq:rw}) - (\ref{eq:MMMaaaaa}), we see
        \begin{eqnarray}\label{eq:r1}
&&\int\limits_{\mathbb{T}}\left|\hat{a}(\theta)\right|^{4}d\theta \ll N^{2} (\log\,N)^{4}\sum\limits_{0<|k|\leq N} d(k)
\sum\limits_{ x,y\leq \sqrt{WN+b},\,x,y\,prime\atop {x^{2}\equiv y^{2}\equiv b(mod\,W)\atop  x^{2}-y^{2}=Wk}}1+N^{3+O(1/\log\log\,N)}.\nonumber\\
           \end{eqnarray} Recalling the standard estimate for the divisor function $d(k)\ll N^{O(1/\log\log N)}$ for $0\leq |k|\leq N,$ we see from (\ref{eq:r1}):
     \begin{eqnarray}\label{eq:r11}
&&\int\limits_{\mathbb{T}}\left|\hat{a}(\theta)\right|^{4}d\theta \ll N^{O(1/\log\log N)}N^{2}(\log\,N)^{4}\sum\limits_{0<|k|\leq N}
\sum\limits_{x,y\leq \sqrt{WN+b},\,x,y\,prime\atop {x^{2}\equiv y^{2}\equiv b(mod\,W)\atop  x^{2}-y^{2}=Wk}}1\nonumber\\
&+&N^{3+O(1/\log\log N)}\nonumber\\
&\ll& N^{2+O(1/\log\log N)}\left(\sum\limits_{x\leq \sqrt{WN+b},\,x\,prime\atop x^{2}\equiv  b(mod\,W)} 1\right)^{2}+N^{3+O(1/\log\log N)}\nonumber\\
&\ll & N^{3+O(1/\log\log N)},
           \end{eqnarray}
which proves Lemma \ref{le:r2}.
\begin{lemma}\label{le:r3}
Define the region
\begin{eqnarray*}
\mathbb{R}_{\delta_{1}}=\left\{\theta\in \mathbb{T}:\left|\hat{a}(\theta)\right|>\delta_{1} N\right\}.
\end{eqnarray*} for any $\delta_{1}\in (0,1).$
For any $\delta_{1}\in (0,1)$ and any $\epsilon>0$ there exists a constant $C_{\epsilon}$ depending only on $\epsilon$ such that
        \begin{eqnarray*}
meas(\mathbb{R}_{\delta_{1}})\leq\frac{C_{\epsilon}}{\delta_{1}^{4+\epsilon}N}.
           \end{eqnarray*}
\end{lemma}
As shown in \cite{browning}, Lemma \ref{le:r1} is a direct consequence of Lemma \ref{le:r3}.
\\
{\it Proof of Lemma \ref{le:r3}:}
Using Lemma \ref{le:r2} instead of \cite[Lemma 6.2]{browning}, we derive in the same way as in \cite[Proof Lemma 6.3]{browning} that we only need to consider the case
\begin{eqnarray}\label{eq:r2}
\delta_{1} >N^{-C(\epsilon \log\log N)^{-1}}.
\end{eqnarray}\\
As in \cite[Section 6]{browning}, we let $\theta_{1},..,\theta_{R}$ be $1/N$ spaced points in $\mathbb{T}$ such that $|\hat{a}(\theta_{r})|\geq \delta_{1} N$, for $1\leq r\leq R.$ We known from \cite[Proof Lemma 6.3]{browning} that in order to prove Lemma \ref{le:r3}, it is sufficient to show that
\begin{eqnarray}\label{eq:hh1}
R\ll \frac{C_{\epsilon}}{\delta_{1}^{4+\epsilon}} .
\end{eqnarray}
In order to prove (\ref{eq:hh1}), we let $f_{n}\in \mathbb{R}$ be such that $|f_{n}|\leq 1$ and $a(n)=f_{n}v(n)$ for integers $1\leq n\leq N.$ Furthemore, define $c_{r}\in \mathbb{C}$ with $|c_{r}|=1$ such that $c_{r} \hat{a}(\theta_{r})=\left|\hat{a}(\theta_{r}\right|$ for $1\leq r\leq R.$  Then it follows from the Cauchy-Schwarz inequality and the prime number theorem in arithmetic progressions with the constant module $W,$
\begin{eqnarray*}
\delta_{1}^{2} N^{2}R^{2} &\leq &\left(\sum\limits_{1\leq r\leq R}|\hat{a}(\theta_{r})|\right)^{2}\nonumber\\
&=& \left(\sum\limits_{1\leq r\leq R}c_{r}\sum\limits_{n\leq N} a(n)e(n\theta_{r})\right)^{2}\nonumber\\
&=& \left(\sum\limits_{n\leq N} a(n)^{1/2}(f_{n}v(n))^{1/2}\sum\limits_{1\leq r\leq R}c_{r}e(n\theta_{r})\right)^{2}\nonumber\\
&\ll & \sum\limits_{n\leq N}a(n)\sum\limits_{n\leq N}v(n)\left|\sum\limits_{1\leq r\leq R}c_{r}e(n\theta_{r})\right|^{2}\nonumber\\
&\ll & N\sum\limits_{n\leq N}v(n)\left|\sum\limits_{1\leq r\leq R}c_{r}e(n\theta_{r})\right|^{2},
\end{eqnarray*}
which implies
\begin{eqnarray}\label{eq:hu}
\delta_{1}^{2}NR^{2}\ll \sum\limits_{1\leq r,r'\leq R}\left|\hat{v}(\theta_{r}-\theta_{r'})\right|.
\end{eqnarray}
Assume $\gamma >2$ be fixed. Applying H\"older's inequality to (\ref{eq:hu}), we find
\begin{eqnarray}\label{eq:kh}
\delta_{1}^{2\gamma}N^{\gamma}R^{2}\ll \sum\limits_{1\leq r,r'\leq R}\left|\hat{v}(\theta_{r}-\theta_{r'})\right|^{\gamma}.
\end{eqnarray}
We put $\theta=\theta_{r}-\theta_{r'}$ for given $r\ne r'.$
With $\theta$ in place of $\frac{r}{N},$ we see from  (\ref{eq:defv1v2}), (\ref{eq:mi23aa}), (\ref{eq:M2u}), (\ref{eq:M31uaaa}), (\ref{eq:M3u}), and Lemma \ref{le:pseudo2},
\begin{eqnarray}\label{eq:m999}
\left|\hat{v}\left(\theta\right)\right|
&\ll &\frac{\phi(W)}{W}\log\,z\left(|T(q_{1}q_{2})|+q^{-1/2+\epsilon}\right)\left|\int\limits_{0}^{N}e\left(\left((\theta-\frac{a}{q}\right) x\right) dx\right|+N^{1-\delta/400}g.\nonumber\\ \end{eqnarray} Using (\ref{eq:kjjj}), Lemma \ref{le:a3}, (\ref{eq:q1}), and (\ref{eq:M31uab}), we see from (\ref{eq:m999}),
\begin{eqnarray}\label{eq:s4}
\left|\hat{v}\left(\theta\right)\right|&\ll &\left(\log\,N|T(q_{1}q_{2})|+q^{-1/2+\epsilon}\right)\min\left\{N,\left|\left|\theta-\frac{a}{q}\right|\right|^{-1}\right\}+N^{1-\delta/400}\nonumber\\
&\ll &N\left(q^{-1/2+\epsilon}\left(1+N\left|\left|\theta-\frac{a}{q}\right|\right|\right)^{-1}+N^{-\delta/400}\right).\nonumber\\
\end{eqnarray}
Let $Q_{1}=\delta_{1}^{-d}$ for some $d\in \left[5,6\right]$ which we will fix later.  If $q> Q_{1},$ then by (\ref{eq:r2}) and (\ref{eq:s4}),
\begin{eqnarray}\label{eq:mumuaaa}\left|\hat{v}\left(\theta\right)\right|\ll Q_{1}^{-1/2+\epsilon}N\leq \delta_{1}^{(5/2+5\epsilon)}N.
\end{eqnarray}
We derive from (\ref{eq:s4})  and (\ref{eq:mumuaaa}) that the contribution of all $\theta=\theta_{r}-\theta_{r'}\in \mathfrak{M}(a,q)$ with $q\geq Q_{1}$ in the sum over $r,\,r'$ in the right-hand side of (\ref{eq:kh}) is $\ll R^{2}\delta_{1}^{\gamma(5/2+5\epsilon)}N^{\gamma},$  which implies that it is smaller than the left-hand side of (\ref{eq:kh}) and we can therefore neglect it.
Summarizing the above, we derive from (\ref{eq:r2}), (\ref{eq:kh}), and (\ref{eq:s4}),
\begin{eqnarray*}
\delta_{1}^{2\gamma}R^{2}& \ll &\sum\limits_{q\leq Q_{1}}\sum\limits_{a\,mod\,q\atop (a,q)=1}\sum\limits_{1\leq r,r'\leq R}\frac{q^{-\gamma/2+\epsilon}}{\left(1+N\left|\left|\theta_{r}-\theta_{r'}-\frac{a}{q}\right|\right|\right)^{\gamma}},
\end{eqnarray*}
which implies
\begin{eqnarray}\label{eq:u88}
\delta_{1}^{2\gamma+\epsilon d}R^{2}&\leq & \sum\limits_{q\leq Q_{1}}\sum\limits_{a\,mod\,q\atop (a,q)=1}\sum\limits_{1\leq r,r'\leq R}\frac{q^{-\gamma/2}}{\left(1+N\left|\left|\theta_{r}-\theta_{r'}-\frac{a}{q}\right|\right|\right)^{\gamma}}.
\end{eqnarray}
We assume $\epsilon\leq 10^{-9}.$ We set $\delta_{2}=\delta_{1}^{\frac{2\gamma+d\epsilon}{2\gamma}}.$ Further we choose $d$ such that
$5=d\frac{2\gamma}{2\gamma+d\epsilon}$ which implies that $d=\frac{10\gamma}{2\gamma-5\epsilon}\in \left[5,6\right].
$ Thus, we see from (\ref{eq:u88}):
\begin{eqnarray}\label{eq:u89}
\delta_{2}^{2\gamma}R^{2}&\leq & \sum\limits_{q\leq \delta_{2}^{-5}}\sum\limits_{a\,mod\,q\atop (a,q)=1}\sum\limits_{1\leq r,r'\leq R}\frac{q^{-\gamma/2}}{\left(1+N\left|\left|\theta_{r}-\theta_{r'}-\frac{a}{q}\right|\right|\right)^{\gamma}}.
\end{eqnarray}
This is exactly the same expression that was found in the proof of \cite[Lemma 6.3]{browning}. It is stated in \cite[end of section 6]{browning} that, using the argument in \cite{bour}, (\ref{eq:u89}) implies (\ref{eq:hh1}) with $\delta_{2}$ instead of $\delta_{1}.$
As $\epsilon \leq 10^{-9}$ and $\gamma>2,$ this in turn implies
\begin{eqnarray*}
R\ll \frac{C_{\epsilon}}{\delta_{1}^{(4+\epsilon)\frac{2\gamma+d\epsilon}{2\gamma}} }\leq \frac{C_{\epsilon}}{\delta_{1}^{4+8\epsilon}}.
\end{eqnarray*} This proves (\ref{eq:hh1}) with $\epsilon:=8\epsilon.$

\section{Regularity condition}\label{se:regularity}
\setcounter{equation}{0}\setcounter{theorem}{0}\setcounter{lemma}{0}
In this section, we prove that the function $a_{4}(n)$ satisfies the regularity condition (5) of Theorem \ref{th:1}. Following the argument in \cite[Section 5]{shao}, we write $\beta=\delta/50,$ $Y=\prod\limits_{p\leq \beta^{-1}\atop p\,is\,prime}p,$ and $\beta $ is chosen such that $Y|W.$ We set
\begin{eqnarray*}
U&=& \left\{1\leq u\leq \beta N_{4}:\,Wu+b_{4}\in \mathbb{P}^{2}\right\},\nonumber\\
V&=& \left\{(1-\beta)N_{4}\leq v\leq  N_{4}:\,Wv+b_{4}\in \mathbb{P}^{2} \right\}.
\end{eqnarray*}
We see,
\begin{eqnarray}\label{eq:k1}
&&\sum\limits_{(u,v)\in M}a_{4}(u)a_{4}(v)=\left(\frac{2\phi(W)}{\sigma(b)W}\log\,z_{4}\right)^{2}\sum\limits_{(u,v)\in M}\sqrt{Wu+b}\sqrt{Wv+b}\nonumber\\
&\gg &NW\left(\frac{\phi(W)}{\sigma(b)W}\log\,z_{4}\right)^{2} \sum\limits_{u\in U,v\in V\atop (u-v,Y)=1}1\nonumber\\
&\gg &NW\left(\frac{\phi(W)}{\sigma(b)W}\log\,z_{4}\right)^{2}\sum\limits_{s_{1},s_{2}(mod\,Y)\atop (s_{1}-s_{2},Y)=1}\left|U\cap Y\mathbb{Z}+s_{1}\right|\left|V\cap Y\mathbb{Z}+s_{2}\right|.
\end{eqnarray}
We now consider the set $U\cap Y\mathbb{Z}+s_{1}.$ For any $u\in U,$ there exists a unique pair of integers $h_{j}$ and $u_{1},$ $h_{j}\in [W],\,j\in [1,..,\sigma(b_{4})],\,h_{j}^{2}\equiv b_{4}(mod\,W),$ $u_{1}\in [1, U_{1}],\,U_{1}:=\left(\sqrt{W\beta N_{4}+b_{4}}-h_{j}\right)/W,$  such that $Wu_{1}+h_{j}\in \mathbb{P}$ and
 \begin{eqnarray}\label{eq:j8j8}\left(Wu_{1}+h_{j}\right)^{2}=Wu+b_{4}.\end{eqnarray} Conversely, for each such pair $u_{1}$ and $h_{j}$ there exists exactly one $u\in U$ such that (\ref{eq:j8j8}) holds.
We write $h_{j}^{2}=c_{j}W+b_{4}$ such that $Wu+b_{4}=W(Wu_{1}^{2}+2h_{j}u_{1}+c_{j})+b_{4}.$ As $Y|W,$ the congruence condition $u\equiv s_{1}(mod\,Y)$ in (\ref{eq:k1}) implies that
\begin{eqnarray}\label{eq:k2}
2h_{j}u_{1}\equiv s_{1}-c_{j}(mod\,Y).
\end{eqnarray}
In the same way, considering the set $V\cap Y\mathbb{Z}+s_{2},$ we find integers $v_{1}\in [V_{1},V_{2}],$ $V_{1}=\left(\sqrt{W(1-\beta)N+b_{4}}-h_{k}\right)/W,$ $V_{2}=\left(\sqrt{WN+b_{4}}-h_{k}\right)/W,$ $h_{k},$ and $c_{k}$ such that $h_{k}^{2}\equiv b_{4}(mod\,W),$ and
\begin{eqnarray}\label{eq:k3}
2h_{k}v_{1}\equiv s_{2}-c_{k}(mod\,Y).
\end{eqnarray}
Combining (\ref{eq:k1}) - (\ref{eq:k3}), we obtain
\begin{eqnarray}\label{eq:k4}
&&\sum\limits_{(u,v)\in M}a_{4}(u)a_{4}(v)\nonumber\\
&\gg &NW\left(\frac{\phi(W)}{\sigma(b)W}\log\,z_{4}\right)^{2}
\sum\limits_{s_{1},s_{2}(mod\,Y)\atop (s_{1}-s_{2},Y)=1}\sum\limits_{j=1}^{\sigma(b)}\sum\limits_{h=1}^{\sigma(b)}\sum\limits_{1\leq u_{1}\leq U_{1}\atop {Wu_{1}+h_{j}\in \mathbb{P}\atop 2h_{j}u_{1}\equiv s_{1}-c_{j}(mod\,Y)}
}\sum\limits_{V_{1}\leq v_{1}\leq V_{2}\atop {Wv_{1}+h_{k}\in \mathbb{P}\atop 2h_{k}v_{1}\equiv s_{2}-c_{k}(mod\,Y)}
}1.\nonumber\\
\end{eqnarray}
For a fixed pair $s_{1}$ and $s_{2}$ with $(s_{1}-s_{2},Y)=1,$ note that the summation condition $(s_{1}-s_{2},Y)=1$ implies that $(s_{1},2)=1$ and $(s_{2},2)=2$ or viceversa. Assuming - without loss of generality - that $(s_{1},2)=1$ and $(s_{2},2)=2$, the congruence conditions in the summations over
$u_{1}$ and $v_{1}$ imply that $(c_{j},2)=1$ and $(c_{k},2)=2.$ Assuming that for each pair $(s_{1},s_{2})$ we can find at least one pair $h_{j0}$ and $h_{k0}$ such that $c_{j0}$ and $c_{k0}$ satisfy these conditions, we see from (\ref{eq:k4}):
\begin{eqnarray*}\label{eq:k5}
&&\sum\limits_{(u,v)\in M}a_{4}(u)a_{4}(v)\nonumber\\
&\geq &NW\left(\frac{\phi(W)}{W}\log\,z_{4}\right)^{2}
\sum\limits_{s_{1},s_{2}(mod\,Y)\atop (s_{1}-s_{2},Y)=1}\sum\limits_{1\leq u_{1}\leq U_{1}\atop {Wu_{1}+a_{j0}\in \mathbb{P}\atop u_{1}\equiv \bar{h_{j0}}(s_{1}-c_{j0})/2(mod\,Y)}
}\sum\limits_{V_{1}\leq v_{1}\leq V_{2}\atop {Wv_{1}+h_{k0}\in \mathbb{P}\atop v_{1}\equiv \bar{h_{k0}}(s_{2}-c_{k0})/2(mod\,Y)}
}1\nonumber\\
&\gg & NW\left(\frac{\phi(W)}{W}\log\,z_{4}\right)^{2}Y\phi(Y)\left(\frac{ \sqrt{\beta N_{4}}}{\log\,N_{4}}\frac{\sqrt{W}}{\phi(W)Y}\right)^{2}\nonumber\\
&\gg & N^{2}\kappa,
\end{eqnarray*}
Here, we have used the prime number theorem in arithmetic progressions of modulus $WY$ which is a constant. To conclude the proof, we show that for a fixed pair $s_{1}$ and $s_{2}$
we can indeed find at least one pair $h_{j0}$ and $h_{k0}$ with $h_{j0}^{2}\equiv h_{k0}^{2}\equiv b_{4}(mod\,W)$ such that $(c_{j0},2)=1$ and $(c_{k0},2)=2.$ Assume a given $h_{j0}$ with $h_{j0}^{2}=c_{j0}W+ b_{4}$ such that $(c_{j0},2)=2.$ We define $h_{j1}:=h_{j0}+\frac{W}{2}.$
We see $h_{j1}^{2}=Wc_{j1}+b_{4},$ where $c_{j1}=c_{jo}+ h_{jo}+\frac{W}{4}.$ As $(c_{jo}+ h_{jo}+\frac{W}{4},2)=1,$ we have found $h_{j1}$ such that $(c_{j1},2)=1.$ In the same way, we can show that there exists $h_{ko}$ such that $(c_{k0},2)=1.$

\end{document}